\documentclass[11pt,reqno]{amsart}
\usepackage[T1]{fontenc}
\usepackage{lmodern}
\usepackage{amsmath,amssymb,amsthm}
\usepackage[letterpaper,textwidth=5.85in,textheight=8.85in,centering]{geometry}
\usepackage[hidelinks,bookmarksopen=false,bookmarksnumbered=true]{hyperref}
\hypersetup{pdftitle={Picard-rank bounds for Calabi--Yau threefolds and the Integral Anosov Conjecture},pdfsubject={Birational geometry of Calabi--Yau threefolds},pdfauthor={Kai Xu}}
\newtheorem{theorem}{Theorem}[section]
\newtheorem{proposition}[theorem]{Proposition}
\newtheorem{lemma}[theorem]{Lemma}
\newtheorem{corollary}[theorem]{Corollary}
\theoremstyle{definition}

\newtheorem{conjecture}{Conjecture}

\theoremstyle{remark}
\newtheorem{remark}[theorem]{Remark}
\numberwithin{equation}{section}
\newcommand{\R}{\mathbb R}
\newcommand{\C}{\mathbb C}
\newcommand{\Mov}{\operatorname{Mov}}
\newcommand{\Eff}{\operatorname{Eff}}
\newcommand{\BigCone}{\operatorname{Big}}

\AddToHook{env/thebibliography/end}{\enlargethispage{2\baselineskip}}
\title[Picard-rank bounds for Calabi--Yau threefolds]{Picard-rank bounds for Calabi--Yau threefolds and the Integral Anosov Conjecture}
\author{Kai Xu}
\address{Yau Mathematical Sciences Center, Tsinghua University, Beijing, China}
\date{20 November 2026, revised 23 November 2026}
\begin{document}
\begin{abstract}
We formulate the Integral Anosov Conjecture and use it to deduce a
uniform bound on the Picard rank of smooth projective strict Calabi--Yau
threefolds whose numerical birational groups act cocompactly on the
whole projectivized movable interior. Together with the finite-chamber
and fibered bounds, this gives a global bound provided that, in
sufficiently large rank, every nonfibered birational class with
infinitely many marked chambers has such an action.
The finite-chamber bound is unconditional. Its proof uses the polytope
of simultaneous divisorial contractions and bounds the independent
choices of exceptional divisors by counting curve classes over their
centres.
\end{abstract}
\maketitle
\raggedbottom
\section{Introduction}

We ask whether there is an absolute upper bound on the Picard number
\(\rho(X)\) of a strict Calabi--Yau threefold. In a smooth M-theory
compactification the number of vector multiplets is \(\rho(X)-1\).
A flop preserves this number and can change the nef cone.

Throughout, a \emph{strict Calabi--Yau threefold} is a smooth projective
complex threefold \(X\) with \(K_X\simeq\mathcal O_X\) and
\(H^1(X,\mathcal O_X)=H^2(X,\mathcal O_X)=0\). Strict transform
identifies the numerical divisor spaces of its smooth projective flop
models with \(N^1(X)_\R\). A \emph{marked chamber} is the nef cone
of such a model in this common space. We include the original model.
A \emph{fibration} is a surjective morphism with connected fibres onto
a normal projective curve or surface.

Write \(\Lambda_X\) for the lattice of integral numerical divisor classes
and \(M=\overline{\Mov}^{\,1}(X)\subset N^1(X)_{\R}\) for the closed
cone generated by classes of divisors with no fixed component. Put
\[
 n=\operatorname{rank}\Lambda_X=\rho(X),\qquad
 \Omega_X=\mathbb P(\operatorname{int}M),\qquad
 \Gamma_X=\operatorname{im}\bigl(\operatorname{Bir}(X)\to
                               \mathrm{GL}(\Lambda_X)\bigr).
\]
We say that \(\Gamma_X\) \emph{divides} \(\Omega_X\) if it acts properly
discontinuously and cocompactly on all of \(\Omega_X\).

Canary--Tsouvalas give the upper bound \(n-1\) for the cohomological
dimension of a torsion-free hyperbolic group with a projective Anosov
inclusion in \(\mathrm{SL}(n,\mathbb R)\), for \(n\geq4\)
\cite[Theorem~1.3]{CanaryTsouvalas}. Our conjecture excludes equality
for integral matrices in all sufficiently large ranks.

\begin{conjecture}[Integral Anosov Conjecture]
\label{conj:integral-anosov}
We conjecture that there is an absolute integer \(N\) with the following
property. Let \(n>N\), and let \(\Gamma\leq\mathrm{SL}(n,\mathbb Z)\)
be any finitely generated, torsion-free, word-hyperbolic subgroup. If
the given inclusion \(\Gamma\hookrightarrow\mathrm{SL}(n,\mathbb R)\)
is projective Anosov, then
\[
                     \operatorname{cd}_{\mathbb Z}(\Gamma)\leq n-2 .
\]
\end{conjecture}

Here \(\operatorname{cd}_{\mathbb Z}\) denotes integral cohomological
dimension, and projective Anosov means \(P_1\)-Anosov. The same \(N\)
must work for every such inclusion. The conjecture requires no
assumption of Coxeter type, commensurability with a reflection group,
an invariant quadratic form, a specified wall structure, or a
Calabi--Yau realization.

\begin{theorem}[Main theorem on boundedness of Picard rank]\label{thm:picard-boundedness}
Assume the Integral Anosov Conjecture with cutoff \(N\).
\begin{enumerate}
\item If \(\Gamma_X\) divides \(\Omega_X\), then
\(\rho(X)\leq\max\{9,N\}\). If \(X\) has finitely many marked
chambers or some smooth projective flop model of \(X\) admits a
fibration, then \(\rho(X)\leq B_0\) for an absolute constant
\(B_0\geq9\) independent of \(N\).
\item Suppose there is an absolute integer \(b\) with the following
geometric property. Whenever \(\rho(X)>b\), \(X\) has infinitely many
marked chambers, and no smooth projective flop model admits a fibration,
\(\Gamma_X\) divides \(\Omega_X\). Then \emph{every} smooth projective
strict Calabi--Yau threefold satisfies
\[
                         \rho(X)\leq\max\{b,B_0,N\}.
\]
\end{enumerate}
\end{theorem}

\begin{proof}
First suppose \(\Gamma_X\) divides \(\Omega_X\).
If \(c_2(X)\) is numerically zero, there is a finite \'etale cover by an
abelian threefold \(A\) \cite{COguisoSakurai}. Pullback is injective
on numerical divisor classes by the projection formula, so
\(\rho(X)\leq h^{1,1}(A)=9\).

Now suppose \(c_2(X)\) is numerically nonzero, and assume for
contradiction that \(n=\rho(X)>\max\{9,N\}\).
Proposition~\ref{prop:c-group} gives strict convexity of \(\Omega_X\).
Lemma~\ref{lem:divided-dimension} then gives a finitely generated,
torsion-free, word-hyperbolic subgroup
\(\Gamma_0\leq\mathrm{SL}(\Lambda_X)\) whose inclusion is projective
Anosov. Its cohomological dimension is the dimension of its closed
aspherical quotient of \(\Omega_X\). Applying our conjecture to this
same inclusion yields
\[
             n-1=\operatorname{cd}_{\mathbb Z}(\Gamma_0)\leq n-2,
\]
a contradiction.

For a fibered flop model \(Y\), let \(F\) be the absolute Picard-rank
bound supplied by boundedness of fibered Calabi--Yau threefolds
\cite[Corollary~1.2]{CFibered}. Then \(\rho(X)=\rho(Y)\leq F\).
For a threefold with finitely many marked chambers and no fibered flop
model, Theorem~\ref{thm:main} below gives the absolute bound \(R\).
Thus \(B_0=\max\{F,R,9\}\) works in (1). Under the hypothesis of (2),
every \(X\) of rank greater than \(b\) satisfies one of the three
conditions of (1), proving the asserted global bound.
\end{proof}

The geometric hypothesis in (2) is the remaining step in extending
this argument to all threefolds. The following consequence gives a case
in which one can verify the required compactness directly.

\begin{corollary}[A bound from the effective movable cone]\label{cor:cone-rank}
Assume the Integral Anosov Conjecture. Suppose no smooth projective
flop model of \(X\) admits a divisorial contraction and the effective
movable cone is covered by the \(\Gamma_X\)-translates of a rational
polyhedral cone contained in it. Then \(\rho(X)\leq\max\{B_0,N\}\).
\end{corollary}
\begin{proof}
If a flop model is fibered, Theorem~\ref{thm:picard-boundedness}(1)
gives \(\rho(X)\leq B_0\). Otherwise Proposition~\ref{prop:division-criterion}
shows that \(\Gamma_X\) divides \(\Omega_X\), so the same theorem gives
\(\rho(X)\leq\max\{9,N\}\).
\end{proof}

\subsection*{The unconditional finite-chamber bound}

\begin{theorem}[Finite-chamber bound]\label{thm:main}
There is an absolute constant \(R\) such that \(\rho(X)\le R\) whenever
\(X\) has finitely many marked chambers and no smooth projective flop
model of \(X\) admits a fibration.
\end{theorem}

The constant \(R\) is independent of the number of marked chambers.
Isomorphic models can have different markings.

Nikulin and Shokurov showed that a single model with \(\rho(X)>40\)
has a small extremal ray or a nonzero real nef class of cube zero
\cite[Theorem~0.1]{NikulinCY}. A small ray can lead to another flop.
Let \(P=\overline{\Eff}^{\,1}(X)\) be the
pseudoeffective cone. Its interior consists of the big classes.

\begin{theorem}[Bound from bigness of movable classes]\label{thm:mov}
The Picard numbers of strict Calabi--Yau threefolds satisfying
\(M\setminus\{0\}\subset\operatorname{int}P\) are uniformly bounded.
\end{theorem}

To deduce Theorem~\ref{thm:main}, we must turn a real boundary class
into a morphism.
With finitely many chambers, a nonzero class on \(M\cap\partial P\)
is nef on some model, but its ray can be irrational. In rank at least
\(19\), Section~\ref{sec:reduction} uses the intersection cubic and
the divisorial walls to obtain a nonzero \emph{rational} nef class of
cube zero with an effective multiple on a smooth flop model.
Threefold log abundance then gives a fibration. The threshold \(19\)
comes from weak approximation on a cubic of dimension \(\rho(X)-2\).

\subsection*{Which divisors can be contracted together?}
Suppose now that every nonzero movable class is big. For a birational
contraction \(f\) on a marked model, put \(D=f^*H\) with \(H\) ample.
The projection formula detects the exceptional primes of \(f\) through
the equality \(D^2 E=H^2 f_*E\). Thus
\[
                D^2E=0\quad\Longleftrightarrow\quad\dim f(E)\le1.
\]
Across all the models, the positive square \(\langle D^2\rangle\)
gives a simple polytope \(Q\) of dimension \(\rho(X)-1\), obtained by
slicing \(P^\vee\). Its facets carry prime divisors as labels. A
collection of facets meets exactly when their primes can be contracted
together. We call such a collection \emph{compatible}.

A simplex and a cube illustrate the two parts of the face argument.
In an \(n\)-simplex any \(n\) facets meet, whereas all \(n+1\) have
empty intersection. For \(Q\), the geometry of the first moving
linear system on a minimally incompatible set bounds its size by ten
in large rank. In an \(n\)-cube the only minimal obstructions are
pairs of opposite facets, and one may choose a side independently in
each direction. Section~\ref{sec:diagram} shows that large dimension
forces many independent choices in \(Q\), in the following precise
sense. There are at most eight common primes \(K\) and disjoint
blocks \(P_1,\ldots,P_m\). With \(K\) present, no whole block can be
contracted. Omit any one prime from each block, and all
the remaining primes can be contracted together with \(K\). For every prescribed \(m\), sufficiently large Picard number forces
at least this many blocks.

\subsection*{How can the choices be counted on one model?}
The contractions just described may live on different models. To
compare them, assign positive weights of sum one in each block.
For a two-prime block \(\{E,E'\}\), a weight \(t\) prescribes the ratio
\[
                    \frac{w(E)}{w(E)+w(E')}=t
\]
for a numerical curve class \(w\). A family consisting of a single
class can realize only one such ratio. More generally, realizing
independent ratios requires parameters.

Section~\ref{sec:weights} retains at least \((m-1)/2\) blocks and
finds one contraction and one \(w\) for all the retained blocks.
The class \(w\) is a convex combination of sums of classes sweeping
a fixed set of at most eight exceptional primes. Each block now counts
in one of two ways. On an untouched block \(w\) realizes the prescribed
ratios. Choosing the weights generally forces the affine dimension of
the available classes to be at least the number of untouched blocks.
On any other block a contracted prime has positive \(w\)-degree, so
it meets a curve sweeping one of the fixed primes.

Both counts are local over the centres of the fixed primes. Over a
curve the sweeping class is fixed, and there are at most three neighbours
in the transverse Du Val diagram. Over a point \(p\), write \(V_p\)
for the numerical span of the curves in the fibre. It contains the
available curve directions over \(p\), and its dimension also bounds
the number of exceptional primes whose centres contain \(p\). Theorem~\ref{thm:fibre}
bounds \(\dim V_p\) in terms of the number of primes mapping to \(p\)
itself. To bound that number, choose the weights near the barycentres
of the blocks. Every contracted block prime then has \(w\)-degree
greater than \(1/20\). Section~\ref{sec:conclusion} bounds their
total contribution by intersecting with a principal divisor through
\(p\). This bounds the number of blocks and proves
Theorem~\ref{thm:mov}.

We use the dual marking on \(N_1\) for curve classes. Intersection cubics
and second Chern classes are taken on the indicated smooth model.
An \emph{elementary} contraction has relative Picard number one.

Appendix~\ref{app:project-results} gives further rank bounds and fibration criteria beyond the finite-chamber setting.

\subsection*{Acknowledgments}
I owe my deepest gratitude to my teachers, Shing-Tung Yau and Cumrun Vafa.
Over many years, they guided my education in mathematics and physics and
led me to the problem studied here. Their example continues to shape the
questions I ask and the way I understand the relationship between geometry
and physics.

I am also grateful to Cumrun Vafa and Hee-Cheol Kim for our very fruitful
collaboration on bounds for the Hodge numbers of elliptically fibered
Calabi--Yau threefolds~\cite{KimVafaXu}.

\section{Exceptional divisors and independent contractions}\label{sec:diagram}

Assume throughout this section that every nonzero class in the closed
movable cone \(M\) is big, and put \(n=\rho(X)-1\ge2\). We encode
simultaneous contractions by the faces of a simple \(n\)-polytope.
The geometric input is that in sufficiently large Picard number, any
failure of simultaneous contractibility is already visible on at most ten primes. The rest of the section
extracts from high dimension many independent choices of primes.

\subsection{The polytope of simultaneous contractions}

Write \(P\) for the closed pseudoeffective cone. Kawamata shows that
\(P\) is generated up to closure by movable classes and exceptional
primes on marked models, and is locally rational polyhedral away
from \(M\) \cite[Theorem~2.9]{Kawamata}. By hypothesis the compact
projectivized boundary of \(P\) lies in this locally polyhedral region.
Thus \(P\) is rational polyhedral, with an exceptional prime on every
extremal ray.

The reason its proper faces are simplicial is rigidity. If the class of
an effective integral divisor \(E\) belongs to a proper face, then
\(h^0(mE)=1\) for every \(m>0\). A nonzero moving part would itself lie
on the face, contrary to the hypothesis on \(M\). Now a dependence among
prime classes on the face would give, after clearing denominators,
numerically equivalent effective divisors with disjoint supports.
Some multiples would be linearly equivalent, since
\(H^1(X,\mathcal O_X)=0\), contradicting rigidity. The primes on the face
are therefore linearly independent. In particular there is one prime
on each extremal ray of the face, and these are all its primes.

Fix an ample integral class \(A\). The dual slice
\[
              Q=P^\vee\cap\{\alpha:\alpha(A)=1\}
\]
is a simple \(n\)-polytope. Each facet bears the label of the prime \(E\)
whose class it annihilates. A set of labels is \emph{compatible} if
their facets meet. To see what compatibility means geometrically, use
the positive-square homeomorphism
\[
       \mathbb P(M)\longrightarrow Q,\qquad
       [D]\longmapsto\frac{\langle D^2\rangle}{\langle D^2\rangle(A)} .
\]
Indeed, \(P^\vee\) is the movable curve cone, and each nonzero element
is strictly positive on \(M\setminus\{0\}\subset\operatorname{int}P\).
The asserted homeomorphism follows from
\cite[Lemma~3.11 and Corollary~3.23]{LehmannXiao}.

Given \(\alpha\in Q\), make the corresponding \(D\) nef on a smooth
marked model. Since \(D\) is big, its nef cone is locally rational
polyhedral there. A rational class in the relative interior of the
same nef face defines a birational contraction \(f\), and \(D=f^*H\)
for a real ample class \(H\) on its target
\cite[Theorems~1.9, 2.3 and 2.6]{Kawamata}. The positive square on
this model is the ordinary square. The projection formula gives
\[
        \alpha(E)=0\quad\Longleftrightarrow\quad
        H^2\cdot f_*E=0\quad\Longleftrightarrow\quad \dim f(E)\le1 .
\]
Thus the labels at \(\alpha\) are exactly the exceptional primes of \(f\).
Conversely, the pullback of an ample class under any such contraction
gives a point with exactly its exceptional labels. In a simple polytope
the relative interior of an intersection of facets lies on precisely
those facets. Hence every compatible set occurs as the \emph{entire set
of exceptional primes} of a contraction. We shall also use that only
finitely many marked chambers occur here, by local finiteness in the
big cone and compactness of \(\mathbb P(M)\)
\cite[Theorem~2.6]{Kawamata}.

\subsection{How many primes can obstruct a contraction?}

A \emph{minimal nonface} is an incompatible set \(E_1,\ldots,E_k\)
all of whose proper subsets are compatible. The sum of the \(E_i\)
is big. Otherwise a nonzero supporting functional would vanish on
the sum, hence on every summand. Every effective divisor supported
on a proper subset is rigid in all multiples.

Choose an effective integral divisor \(G=\sum b_i E_i\) with at least
two sections and with \(\sum b_i\) as small as possible. Such a divisor
exists by bigness of the sum. Rigidity on proper subsets gives
\(b_i>0\) for all \(i\), and minimality gives
\begin{equation}\label{eq:minimal-system}
                 h^0(G)\ge2,\qquad h^0(G-E_i)=1\quad(1\le i\le k).
\end{equation}
In particular \(|G|\) has no fixed component, since removing one would
preserve all sections. Thus \(G\) is movable and big. On its ample
model, \eqref{eq:minimal-system} makes the image of \emph{each} prime
carry almost all the sections of \(G\). We shall bound how many such
images the Riemann--Roch formula allows.

Write \(B(v)\) for a bound on the Picard numbers of smooth projective
Calabi--Yau threefolds carrying a big nef integral divisor of cube
at most \(v\) \cite[Corollary~1.6]{Birkar}. The estimate for curve
images comes from a transverse surface section.

\begin{lemma}[The contribution of a curve centre]\label{lem:chern}
Let \(f:Y\to Z\) be the birational contraction defined by a big rational
nef class on a smooth projective Calabi--Yau threefold. For each curve
centre \(C\), let \(r_C\) be the number of exceptional curves in a
general transverse Du Val resolution and \(\Gamma_C\) its local
quotient group. For every ample rational Cartier class \(H\) on \(Z\),
\[
       c_2(Y) f^*H\ \ge\ \sum_C (H\cdot C)\bigl(r_C+1-|\Gamma_C|^{-1}\bigr).
\]
\end{lemma}
\begin{proof}
For sufficiently divisible \(m\), a general member of \(|mH|\) avoids
the point centres and meets each \(C\) in \(m(H\cdot C)\) general points.
Its pullback replaces each of these points, of orbifold Euler number
\(1/|\Gamma_C|\), by a tree of Euler number \(r_C+1\). Adjunction gives
the same normal-bundle contribution on both surfaces, since the normal
bundle upstairs is pulled back. Dividing their Euler difference by
\(m\) therefore gives
\[
    c_2(Y)f^*H-\widehat c_2(Z)H
       =\sum_C(H\cdot C)\bigl(r_C+1-|\Gamma_C|^{-1}\bigr).
\]
Here \(\widehat c_2\) is the orbifold second Chern class.
The crepant target has \(K_Z\sim0\), so
\(\widehat c_2(Z)H\ge0\) by \cite[Theorem~6]{CGG}.
\end{proof}

\begin{lemma}[At most ten divisors in a minimal obstruction]\label{lem:ten-labels}
If \(\rho(X)>B(36)\), every minimal nonface has at most ten labels.
\end{lemma}
\begin{proof}
Make the divisor \(G\) in~\eqref{eq:minimal-system} nef on a smooth
flop model \(X'\). Flops identify the spaces of sections, since the
pole-order conditions defining a section involve only prime divisors.
Writing \(h=h^0(G)\), we find that the restriction of sections to each
\(E_i\) has rank \(h-1\).

Let \(f:X'\to Z\) be the ample model of \(G\). Write \(G=f^*H\) with
\(H\) ample Cartier. Indeed all sufficiently large multiples of \(G\)
are free and descend, so the difference of consecutive such multiples
also descends. Put \(d=G^3\). Since \(\rho(X')=\rho(X)>B(36)\), we have
\(d\ge37\). Riemann--Roch, vanishing and Chern positivity
\cite{Miyaoka} give
\[
          12h=2d+c_2(X')G,\qquad d\le6h,\qquad h\ge7 .
\]
Each image \(W=f(E_i)\) carries at least \(h-1\) independent sections of
\(H|_W\), so it has positive dimension. By nonnegativity of Fujita's
\(\Delta\)-genus \cite[Theorem~1.1]{FujitaDelta},
\[
       H^{\dim W}\cdot W\ge h^0(W,H|_W)-\dim W\ge h-1-\dim W.
\]

Let \(s\) primes have surface image and \(t\) have curve image. In
\(d=\sum b_iG^2E_i\) only surface images contribute, each by at least
\(h-3\). Thus \(1\le s\le d/(h-3)\le10\). Each curve image has
\(H\)-degree at least \(h-2\). If \(t_C\) of the primes map onto \(C\),
the transverse tree has at least \(t_C\) vertices. Since
\(|\Gamma_C|\ge r_C+1\) in every ADE type, Lemma~\ref{lem:chern} gives
\[
      \frac{c_2(X')G}{h-2}
      \ge \sum_{C:\,t_C>0}\left(t_C+\frac{t_C}{t_C+1}\right)
      \ge t+\frac{t}{t+1}.
\]
For the last inequality use \(t_C+1\le t+1\) and \(\sum_C t_C=t\).

If \(s+t\ge11\), substitute \(t\ge11-s\) into this inequality and use
Riemann--Roch. The resulting upper bound and the surface lower bound are
\[
       \frac{2d-24}{h-2}\le s+\frac1{12-s},
       \qquad s(h-2)\le d+s.
\]
Together they imply
\[
       d\le 24+s+\frac{24+2s}{s(12-s)-1}.
\]
For \(1\le s\le10\) we have
\(s(12-s)-1=(s+9)+(s-1)(10-s)\ge s+9\). The last fraction is
less than three, so the integer \(d\) is at most \(26+s\le36\),
a contradiction.
\end{proof}

\subsection{From small obstructions to cubical faces}

Suppose \(\rho(X)>B(36)\) and put \(q=10\). Lemma~\ref{lem:ten-labels}
says that a collection of facets meets as soon as every subcollection
of at most \(q\) meets. For example, in a simplex face, omitting any one
side leaves the opposite vertex. A minimal nonface among its sides
and its containing facets must therefore use every side. Such a simplex
has dimension at most \(q-1\).

The same incidence criterion bounds the number of triangles through an
edge. A count of the sides of all polygons will then show that almost
every face of fixed dimension is a cube. Throughout, a cube means a
polytope combinatorially equivalent to a cube.

\begin{lemma}[Cubical faces in high dimension]\label{lem:large-cube}
For each fixed \(t\ge2\), the polytope \(Q\) has a cubical \(t\)-face
once \(n\) is sufficiently large in terms of \(t\).
For \(t\ge q\), it suffices that \(n>C_qt^{q-1}\) for a constant \(C_q\).
\end{lemma}
\begin{proof}
Write the facets at the endpoints of an edge as \(A\cup\{a\}\) and
\(A\cup\{b\}\). Their union is incompatible. A minimal nonface in it
uses \(a,b\) and at most \(q-2\) members of \(A\). If the edge is in a
triangle, the third vertex lies on \(a,b\) and every member of \(A\)
except one. The omitted member must belong to the minimal nonface.
It determines the two-face containing the triangle, so there are at
most \(q-2\) such triangles through the edge.

Let \(V\) be the vertex count and \(f_j\) the number of \(j\)-faces.
For a generic linear height, a vertex with \(e\) descending edges is
the maximum of exactly \(\binom ej\) faces of dimension \(j\). Average
this count for the height and its negative. Convexity of
\(e\mapsto\binom ej\) gives Nikulin's estimate
\cite[Lemma~1.3]{Nikulin}
\[
       f_j\ge\frac V2\left(\binom{\lfloor n/2\rfloor}j+
                              \binom{\lceil n/2\rceil}j\right).
\]
Each vertex belongs to \(\binom n2\) polygons. The estimate at \(j=2\)
therefore gives the following bound on their total excess over four sides.
\[
       \sum_{\dim F=2}(\#\operatorname{vert} F-4)
                    =V\binom n2-4f_2\le Vn/2 .
\]
Only triangles contribute negatively. The edge bound and the \(nV/2\)
edges of \(Q\) give \(O_q(nV)\) triangles, so at most \(O_q(nV)\)
polygons in all have a number of sides different from four.

By simplicity each polygon lies in \(\binom{n-2}{t-2}\) faces of
dimension \(t\). The number of \(t\)-faces containing a nonsquare
polygon is thus \(O_{q,t}(Vn^{t-1})\), whereas the height estimate
gives at least \(V\binom{\lfloor n/2\rfloor}{t}\) faces in total.
For large \(n\), some \(t\)-face has only quadrilateral two-faces.
A simple polytope with this property is a cube
\cite[Lemma~4.6]{MasudaPanov}. The uniform estimate in \(t\) follows
in Appendix~\ref{app:cube-rate}.
\end{proof}

\subsection{Independent choices with few common primes}

At a vertex of a cubical face one may independently choose a side in
each direction. All facets containing the face remain present, and
there may be many of them. The following proposition keeps at most
eight common primes by allowing a direction to have more than two
labels.

For example, in a product of triangles take the cube obtained by
choosing one edge in each triangle. Its containing facets record the
chosen edges, one per factor. In the full product all three edge labels
of a triangle form a minimal nonface, while any proper subset meets.
Using these triples as blocks removes every common label and preserves
independence across factors.

\begin{proposition}[Independent proper subsets]\label{prop:blocks}
There is a universal polynomial \(p\) with the following property.
If \(m\ge1\), \(\rho(X)>B(36)\) and \(n>p(m)\), there are pairwise
disjoint sets of labels \(K,P_1,\ldots,P_m\) with \(|K|\le8\) and
\(2\le |P_i|\le10\). Each \(K\cup P_i\) is a minimal nonface, and
\(K\cup\bigcup_i S_i\) is compatible whenever \(S_i\subsetneq P_i\)
for every \(i\).
\end{proposition}
\begin{proof}
Let \(A\) be the facets containing a large cubical face. Denote the
ambient facets cutting out the two opposite sides in direction \(i\)
by \(a_i,b_i\). Among \(A\cup\{a_i,b_i\}\) there is a minimal
nonface. Since either side meets all of \(A\), it has the form
\(B_i\cup\{a_i,b_i\}\), with \(B_i\subseteq A\) and \(|B_i|\le q-2\).

We first arrange for these nonfaces to have disjoint parts outside a
fixed set. The sets \(B_i\cup\{a_i\}\) are distinct and have at most
\(q-1\) labels. By the sunflower lemma \cite{ErdosRado}, at least
\(N\) of them have the same pairwise intersection \(K_0\subseteq A\)
once the cube has \(C_qN^{q-1}\) directions. For those directions put
\[
            D_i=(B_i\setminus K_0)\cup\{a_i,b_i\}.
\]
The \(D_i\) are disjoint, and each \(K_0\cup D_i\) is a minimal nonface.
Moreover, every incompatible set on these labels contains both sides
of some direction. Indeed, any choice using at most one side per
direction lies at a cube vertex together with all the facets in \(A\).

There can still be incompatibilities involving parts of several
\(D_i\). To find independent choices, consider all minimal nonfaces
\(E\) on \(K_0\) and the current \(D_i\). For each \(E\) record the size
\(s\) of a largest part \(E\cap D_i\) and the number \(k\) of parts of
that size. Choose the least pair \((s,k)\) in lexicographic order.
Here \(s\ge2\) by the opposite-side observation. Removing a largest
part of such an \(E\) leaves a set \(K\) with \(|K|\le q-s\le8\).

The counting argument in Appendix~\ref{app:repeated-part} shows that,
after restricting the index family and taking the least pair there,
one can keep the \emph{same} remainder \(K\) and complete it to
a minimal nonface of this least pair in \(m\) distinct sets \(D_i\).
Indeed, if every remainder has few completions, these nonfaces are
sparse enough to eliminate on a smaller family. The least pair then
increases, and only boundedly many pairs are possible. A polynomial
initial \(N\) suffices.

Call the \(m\) alternative parts \(P_1,\ldots,P_m\). Each has size
\(s\), and each \(K\cup P_i\) is a minimal nonface. If we take a proper
subset of every \(P_i\), the union with \(K\) has no part larger than
\(s\), and only the \(k-1\) largest parts in \(K\) can still have size
\(s\). Any minimal nonface in this union would have a smaller pair.
Thus every such choice is compatible. The sunflower bound and
Lemma~\ref{lem:large-cube} give a polynomial threshold for \(n\).
\end{proof}

Equivalently, after fixing \(K\), the compatible subsets of the block
labels form the join of the simplex boundaries \(\partial\Delta(P_i)\).
This is why the next section uses a product of simplices. The zero
coordinates of a point in each factor can be any proper subset of its
labels, and the corresponding primes can all be contracted with \(K\).

\section{Realizing prescribed degree ratios}\label{sec:weights}

Consider first a block of two primes. Give them weights \(t,1-t\).
If both primes survive a contraction, we want a curve class whose
degrees on them are in the ratio \(t:(1-t)\). If one prime is
contracted, we only need its share of the total degree to be at least
its prescribed weight. In the first case the ratio records a parameter.
In the second case the positive degree tells us where to find the
contracted prime. Our aim is to realize this alternative for many
blocks on a single contraction.

\subsection{Sweeping curves and negativity}

For an elementary contraction of a prime \(E_q\), let \(\ell_q\) be the
primitive integral generator of the contracted curve ray, oriented in
the effective direction. We call it a \emph{sweeping class}. The contracted curves cover \(E_q\). A general one has nonnegative
degree on every other prime and on each movable divisor.
Under the marking, \(\ell_q\) therefore acts on \(N^1(X)_\R\) as a
functional nonnegative on \(M\) and positive on its interior. Its sole
negative prime is \(E_q\), with \(1\le-\ell_q(E_q)\le3\) by
Lemma~\ref{lem:degree}.

Suppose a big rational nef class \(D\) defines a contraction \(f:Y\to Z\).
A sweeping class is \emph{available at \(f\)} if its elementary
contraction is realized on a smooth model over \(Z\). In particular
the class annihilates \(D\). Every exceptional prime \(E_q\) supplies at
least one such class. Indeed, run the relative minimal model program
for \((Y,\eta E_q)\) with \(0<\eta\ll1\). The negativity lemma keeps the
strict transform of \(E_q\) from being relatively nef while it survives.
The small steps are flops between smooth models \cite{Kollar,Fujino}.
At the first divisorial step its curves have negative degree on
\(E_q\), so they lie in \(E_q\) and the contracted divisor is \(E_q\).
Termination guarantees that this step is reached.

We will vary \(D\) in \(\mathbb P(M)\). The availability locus of a sweeping
class is the union, over models realizing it, of the closed nef facets
dual to its curve ray. At a point of such a facet,
the contraction defined by \(D\) factors through the elementary one.
Conversely, a model on which the elementary contraction takes place
over \(Z\) is \(D\)-nef, since \(D\) pulls back from \(Z\). The section
ring identifies the ample model across the markings. For a real ray,
use a rational point in the relative interior of the same chamber
face, which lies on exactly the same facets. Section~\ref{sec:diagram}
gives finitely many chambers and facets. Thus the list of sweeping
classes is finite, and the availability locus of each is closed.

\begin{lemma}[Negativity on the swept primes]\label{lem:row-witness}
Let a collection of sweeping classes annihilate a common big class.
Write \(I\) for the primes they sweep. Some divisor
\(V=\sum_{q\in I}v_qE_q\), with all \(v_q>0\), has negative degree on
every class. Moreover, every nonzero effective divisor supported on
\(I\) has negative degree on at least one of the classes.
\end{lemma}
\begin{proof}
Write the big class as \(A+F\), where \(A\) is ample and \(F\) is
effective. Each class has positive degree on \(A\), hence negative
degree on \(F\). Discard from \(F\) all components outside \(I\). Their
degrees are nonnegative, so the divisor \(V\) left over still has
negative degree on every class. In particular, it contains every
swept prime.

Given a nonzero effective \(U=\sum u_q E_q\) supported on \(I\), take
the smallest \(t\) for which \(tV\ge U\). Then \(t>0\) and some
\(E_q\) has coefficient zero in \(tV-U\). Any class sweeping that prime
has nonnegative degree on \(tV-U\), so \(\ell_q(U)\le t\ell_q(V)<0\).
\end{proof}

\subsection{A common set of primes meeting every block}

Start with \(d\ge2\) blocks from Proposition~\ref{prop:blocks}, with
common set \(K\). We want to sum one sweeping class from each common
prime. For the normalized degrees of this sum to make sense, it must
have positive total degree on every block, regardless of which
available classes we choose. The movable divisor on a minimal nonface,
constructed in~\eqref{eq:minimal-system}, will supply this positivity.

If \(K\ne\varnothing\), take \(K_0=K\) and retain all the blocks. On any
contraction of \(K\), choose one available class for each of its primes.
Fix a block \(P_i\) and split the movable divisor supported on
\(K\cup P_i\) into its two parts. By Lemma~\ref{lem:row-witness} some
chosen class is negative on the part supported on \(K\). It is
nonnegative on the whole movable divisor, so it is positive on a
prime of \(P_i\).

If \(K=\varnothing\), each block is itself a minimal nonface. Set aside
one block \(P_0\). We claim that for each other block \(P_i\) there
is a prime \(r_i\in P_0\) such that \emph{every} class sweeping \(r_i\)
is positive on some member of \(P_i\). Otherwise, for every prime
of \(P_0\) we could choose a sweeping class vanishing on \(P_i\).
All would annihilate the big divisor \(\sum_{E\in P_i} E\).
Lemma~\ref{lem:row-witness} would then make one of them negative on
the movable divisor supported on \(P_0\). To obtain a proper subset
of \(P_0\) while retaining many blocks, split \(P_0\) in two parts of
sizes differing by at most one. Take as \(K_0\) the part containing
\(r_i\) for at least half the indices, and retain those \(P_i\).
Since \(K_0\subsetneq P_0\), it is compatible with any simultaneous
choice of proper subsets in the retained blocks.

In both cases the common set \(K_0\), or \emph{core}, has between one
and eight primes. There remain \(b\ge(d-1)/2\) blocks
\(P_i=\{E_{i1},\ldots,E_{im_i}\}\). On any contraction of \(K_0\),
choose one available class per core prime and set
\begin{equation}\label{eq:normal}
                         w=\sum_{q\in K_0}\ell_q .
\end{equation}
The classes may be realized on different marked models over the
same target. All their degrees outside \(K_0\) are nonnegative, and
each block contributes a positive integral degree to at least one
summand. Consequently
\begin{equation}\label{eq:mass}
 w(E_q)\ge-3\quad(q\in K_0),\qquad
 S_i(w):=\sum_j w(E_{ij})\ge1\quad(1\le i\le b).
\end{equation}
The same inequalities hold for convex combinations of these sums.

\subsection{Balancing the degrees on a product of simplices}

Use one simplex per retained block, and put
\(B=\prod_{i=1}^b\Delta^{m_i-1}\). A point \(x\in B\) specifies the
contraction of the common primes and of the primes at its zero coordinates.
Their labels form the set
\[
                  A(x)=K_0\cup\{E_{ij}:x_{ij}=0\}.
\]
We first choose a movable ray \(D(x)\) continuously in \(x\) whose
entire exceptional set is \(A(x)\). To do so, associate to a face \(F\)
of \(B\) the intersection \(Q_F\) of the facets of \(Q\) labelled by
\(A(x)\), for \(x\in\operatorname{relint}F\). Compatibility makes
\(Q_F\) nonempty, and simplicity of \(Q\) makes these its precise
containing facets. The assignment preserves face inclusion. Sending
each barycentre of \(F\) to the barycentre of \(Q_F\) and interpolating
on the barycentric subdivision gives a continuous map \(B\to Q\).
It takes \(\operatorname{relint}F\) into \(\operatorname{relint}Q_F\),
since the barycentre of \(Q_F\) occurs there with positive coefficient.
Composing with the inverse positive-square map gives \(D(x)\).

Let \(C(x)\) be the convex hull of all sums~\eqref{eq:normal} available at
\(D(x)\), and \(W\) the convex hull of all such sums over \(B\).
For \(w\in W\) its normalized degrees on a block form the probability vector
\[
                 u_i(w)=\bigl(w(E_{ij})/S_i(w)\bigr)_j.
\]
By~\eqref{eq:mass} this vector depends continuously on \(w\).
The sets \(C(x)\) are nonempty polytopes in the compact polytope \(W\).
Near any given \(x\), only sums already available at \(x\) can occur,
by finiteness and closed availability. Hence \(x\mapsto C(x)\) has
closed graph.

\begin{proposition}[Simultaneous degree ratios]\label{prop:weights}
For every interior \(y\in B\) there are \(x\in B\), \(w\in C(x)\) and
numbers \(\lambda_i\ge1\) such that
\begin{align}
 x_{ij}=0&\quad\Longrightarrow\quad
      w(E_{ij})=\lambda_i y_{ij}S_i(w),\label{eq:weights}\\
 x_{ij}>0\ \text{for all }j&\quad\Longrightarrow\quad
      w(E_{ij})=y_{ij}S_i(w)\ \text{for all }j.\label{eq:unaffected}
\end{align}
If \(y_{ij}>1/(2m_i)\) for every \(i,j\), each contracted block prime
has degree greater than \(1/20\).
\end{proposition}

\begin{proof}
Fix an interior point \(y_i\) of a simplex. Given a degree vector
\(u_i\ne y_i\), draw the ray from \(u_i\) through \(y_i\) and take its
exit point \(x_i\) on the boundary. The point \(y_i\) is then a strict
convex combination of \(u_i\) and \(x_i\). On every zero coordinate
of \(x_i\), the coordinate of \(u_i\) is consequently the same multiple
of the prescribed coordinate of \(y_i\). If \(u_i=y_i\), the prescribed
degrees are already realized and we allow every \(x_i\) in the simplex.

\begin{center}
\setlength{\unitlength}{0.78pt}
\begin{picture}(170,107)(-7,-16)
  \thinlines
  \put(0,0){\line(1,0){150}}
  \put(0,0){\line(5,6){75}}
  \put(75,90){\line(5,-6){75}}
  \thicklines
  \put(65,52){\vector(3,-4){39}}
  \put(65,52){\circle*{4}}
  \put(83,28){\circle*{4}}
  \put(104,0){\circle*{4}}
  \put(60,53){\makebox(0,0)[r]{\(u_i\)}}
  \put(89,30){\makebox(0,0)[l]{\(y_i\)}}
  \put(110,6){\makebox(0,0)[l]{\(x_i\)}}
  \put(65,-12){\makebox(0,0){\scriptsize \(j\)-th coordinate zero}}
\end{picture}
\end{center}

Write \(R_i(u_i)\) for the singleton exit point in the first case and
the whole simplex in the second. The exit point varies continuously
away from \(y_i\) (radial projection from an interior point of a
polytope). Allowing the whole simplex at \(y_i\) makes the graph of
\(R_i\) closed. On \(B\times W\) consider the correspondence
\[
 (x,w)\longmapsto
       \left(\prod_i R_i\bigl(u_i(w)\bigr)\right)\times C(x).
\]
It has closed graph and nonempty compact convex values. Kakutani's
fixed point theorem \cite{Kakutani} yields \(x_i\in R_i(u_i(w))\)
and \(w\in C(x)\). If \(u_i(w)=y_i\), take \(\lambda_i=1\).
Otherwise write \(y_i=t_i u_i(w)+(1-t_i)x_i\) with \(0<t_i<1\).
On a contracted coordinate, \(u_{ij}(w)=y_{ij}/t_i\), which gives
\eqref{eq:weights} with \(\lambda_i=1/t_i\). An interior \(x_i\)
cannot be an exit point, so it forces \(u_i(w)=y_i\), proving
\eqref{eq:unaffected}. Finally, a contracted prime has degree at least
\(y_{ij}S_i(w)>1/(2m_i)\ge1/20\) for the stated choice of weights.
\end{proof}

\section{Counting blocks over the centres of exceptional divisors}\label{sec:conclusion}

Fix the prescribed weights before taking the contraction. Each block
left untouched by the contraction must realize an independently chosen
ratio. Every other block supplies an exceptional prime meeting a curve
that sweeps a core prime. We count both contributions over the centres
of the core primes.

\subsection{What an untouched block costs}

Choose weights in the open neighbourhood of the barycentre given by
\(y_{ij}>1/(2m_i)\). For a solution \((x,w)\) of
Proposition~\ref{prop:weights} let \(U\) be the indices of the untouched
blocks and put \(L_x=\operatorname{aff} C(x)\). On an untouched block set
\(T_i=\{1,\ldots,m_i\}\). On any other block let \(T_i\) index the
contracted primes. In either case \(T_i\) is nonempty and the ratio
of the degrees on \(T_i\) agrees with the ratio of the prescribed
weights there.

\begin{lemma}[Dimension of the prescribed weights]\label{lem:weight-dimension}
For general rational weights in this neighbourhood, every solution satisfies
\[
                 \sum_i (|T_i|-1)\le\dim L_x .
\]
In particular, \(\sum_{i\in U}(m_i-1)\le\dim L_x\).
\end{lemma}

\begin{proof}
For fixed \((T_i)_i\) and \(L=L_x\), the available ratios are in the
image of the rational map
\[
 L\dashrightarrow\prod_i\mathbb P^{|T_i|-1},\qquad
 w\longmapsto\bigl([w(E_{ij})]_{j\in T_i}\bigr)_i .
\]
If the asserted inequality fails, the closure of this image is a
proper algebraic subset of the target. The prescribed ratios on the
\(T_i\) vary independently, since the positive entries on each selected subset
may be chosen in arbitrary proportions before the block is normalized.
Thus the corresponding map from weight space is dominant, and the
exceptional weights lie in a proper algebraic subset. There are
finitely many choices of \((T_i)_i\) and of \(L\), since the list of
available sums is finite. Choose rational weights in the stated
neighbourhood outside the union of these exceptional sets.
\end{proof}

\subsection{Counting over a fixed contraction}

Choose these weights and a solution. Let \(f:Y\to Z\) be the resulting
contraction. A core prime with curve centre \(\Gamma\) supplies a
single available class. Indeed, a general transverse section gives
the minimal resolution of a Du Val singularity. The fibre components
in the prime form one monodromy orbit and therefore have the same
degrees on marked divisors. Flops over \(Z\) are isomorphisms over a
general point of \(\Gamma\). Thus the primitive generator of this ray
is fixed. Let \(a\) be the number of such core primes and \(w_0\) the
sum of their classes.

Let \(\mathcal P\) be the set of distinct point centres of the remaining
core primes. For \(p\in\mathcal P\), write \(V_p\) for the span of
curves in \(f^{-1}(p)\) and \(v_p=\dim V_p\). A class supplied by a
core prime over \(p\) lies in \(V_p\) on its own model, and
Theorem~\ref{thm:fibre} identifies this space on every marked model
over \(Z\). Hence
\begin{equation}\label{eq:ball}
 C(x)\subseteq w_0+\sum_{p\in\mathcal P} V_p,\qquad
 |U|\le\dim L_x\le\sum_{p\in\mathcal P}v_p .
\end{equation}
Here each untouched block contributes at least one ratio by
Lemma~\ref{lem:weight-dimension}.

For each block outside \(U\), choose one contracted prime \(E\).
It has \(w(E)>1/20\), so some class in the convex combination giving
\(w\) is positive on \(E\). A general curve sweeping its core prime
\(E_q\) therefore meets \(E\), and \(f(E_q)\subseteq f(E)\).
If \(f(E_q)=\Gamma\) is a curve then \(f(E)=\Gamma\). The two primes
give neighbouring orbits in the transverse Du Val diagram. Every
neighbouring orbit has a vertex adjacent to a fixed representative
of the core orbit, by transitivity. Each vertex has valency at most
three. The curve-centred core primes therefore account for at most
\(3a\) chosen primes.

If \(f(E_q)=p\) is a point, then \(E\) maps to \(p\) or to a curve
through \(p\). Write \(k_p\) and \(c_p\) for the total numbers of
exceptional primes of these two types. The chosen primes are distinct,
and the lower bound of Theorem~\ref{thm:fibre} gives
\begin{equation}\label{eq:fibre-use}
             b-|U|\le 3a+\sum_{p\in\mathcal P}(k_p+c_p)
                    \le 3a+\sum_{p\in\mathcal P}v_p .
\end{equation}
Both counts now involve only \(a\) and the dimensions \(v_p\). By the
upper bound of Theorem~\ref{thm:fibre}, it remains to control \(k_p\).

\subsection{The degree of a local equation}

\begin{lemma}[Counting contracted primes over a point]\label{lem:hyperplane}
For \(p\in\mathcal P\), let \(s_p\) be the number of core primes with
centre \(p\). There is an absolute integer \(N\) such that
\(k_p\le(1+60N)s_p\).
\end{lemma}

\begin{proof}
Choose a local equation \(g\) through \(p\) whose vanishing order
\(h_E\) along every prime with centre \(p\) lies between \(1\) and an
absolute integer \(N\). We first explain why such an equation exists.
The volume form on \(Y\) has divisor zero on \(Z\), so \(K_Z\) is
Cartier, \(f\) is crepant and \((Z,p)\) is canonical. A Koll\'ar
component over \(p\) and bounded local complements give a log
canonical pair \((Z,\Delta^+)\) with this component as a log canonical
place and \(N(K_Z+\Delta^+)\sim0\) near \(p\), for \(N\) depending
only on the dimension \cite[Lemma~2.7 and Theorem~2.14]{Xu}.
The support of \(\Delta^+\) passes through \(p\) since \((Z,0)\) is klt.
Take \(\operatorname{div}(g)=N\Delta^+\) locally. The order \(h_E\) is
positive, and log canonicity gives \(h_E\le N\) because every prime on
\(Y\) with centre \(p\) has log discrepancy one over \(Z\).

Keep only the summands of \(w\) sweeping primes over \(p\), with their
original coefficients, and call the result \(w_p\). A general curve
from any other centre misses the primes over \(p\), so \(w_p(E)=w(E)\)
on those primes. On any model realizing a summand, write the pullback as
\[
               \operatorname{div}(g\circ f)
                      =R+\sum_{f(E)=p}h_E E .
\]
The coefficients are unchanged by a flop over \(Z\). The remainder
\(R\) is effective and contains no prime over \(p\), so every summand
of \(w_p\) has nonnegative degree on it. The full pullback has degree
zero on curves over \(p\). Using \(w(E)\ge-3\) for core primes gives
\begin{equation}\label{eq:hyper}
 \sum_{\substack{f(E)=p\\E\notin K_0}}h_Ew(E)
  =-w_p(R)-\sum_{\substack{f(E)=p\\E\in K_0}}h_Ew(E)
  \le3Ns_p .
\end{equation}
Every noncore prime over \(p\) occurs on the left. By the construction
of \(D(x)\), each is indexed by a zero coordinate in a retained block. Each contributes more than \(1/20\) by
Proposition~\ref{prop:weights}. Thus \(k_p-s_p\le60Ns_p\).
\end{proof}

There are at most eight core primes, so \(a+|\mathcal P|\le8\) and
\(k_p\le8(1+60N)\) for each \(p\in\mathcal P\). Put
\(G_*=\max_{1\le k\le8(1+60N)}G(k)\). Adding the two block counts
\eqref{eq:ball} and \eqref{eq:fibre-use} gives
\begin{equation}\label{eq:dimension}
       b\le3a+2\sum_{p\in\mathcal P}v_p
        \le8\max\{3,2G_*\}=:C .
\end{equation}
The dimension of the full product of weight simplices satisfies the
bound \(27a+10\sum_p v_p\), since the untouched blocks contribute at most
\(\dim L_x\), and each remaining block contributes at most nine.

\begin{proof}[Proof of Theorem~\ref{thm:mov}]
Let \(m=2(C+1)+1\) and let \(p\) be the polynomial in
Proposition~\ref{prop:blocks}. If
\(\rho(X)>\max\{B(36),1+p(m)\}\), that proposition supplies \(m\)
blocks. At least \((m-1)/2=C+1\) survive the choice of core, contrary
to~\eqref{eq:dimension}.
\end{proof}

\section{The numerical span of a contracted fibre}\label{sec:local}

Let \(f:Y\to Z\) be the connected-fibre contraction defined by a big
rational nef divisor \(D\) on a smooth projective strict Calabi--Yau
threefold. For \(p\in Z\) a closed point, let \(V_p\subset N_1(Y)_{\R}\)
be the span of the classes of curves in \(f^{-1}(p)\). Let \(k\) be
the number of exceptional primes with centre \(p\), and \(c_p\) the
number whose centre is a curve through \(p\).

\begin{theorem}[Flop invariance and dimension of the fibre curve space]\label{thm:fibre}
The marked space \(V_p\) is independent of the smooth flop model over
\(Z\). Whenever \(k>0\),
\[
                    k+c_p\le\dim V_p\le G(k)
\]
for a function \(G:\mathbb Z_{>0}\to\mathbb Z_{>0}\) independent of
the contraction. One may take \(G(1)=9\), \(G(2)=12\) and \(G(3)=262\).
\end{theorem}

A flop can change the surfaces in the fibre and their individual
intersection numbers. The span of the fibre curves survives the flop,
so we may estimate it on whichever model makes those surfaces easiest
to control. We first relate the contribution of a surface to its
self-intersection. We then choose a model on which all the surface
cubes admit simultaneous lower bounds.

\subsection{What one exceptional surface contributes}

\begin{lemma}[Curve classes on a rigid divisor]\label{lem:degree}
Let \(E\) be a prime divisor rigid in every multiple, so that
\(h^0(Y,\mathcal O_Y(mE))=1\) for every \(m>0\). Put
\(B_E(L,M)=ELM\). The span in \(N_1(Y)_\R\) of curves on \(E\)
has dimension \(\operatorname{rank}B_E\),
and
\[
                    \operatorname{rank} B_E+E^3\le10 .
\]
At an elementary contraction of \(E\), let \(\ell\) be the primitive
integral functional on the contracted ray, oriented effectively.
Then \(1\le-\ell(E)\le3\), and \(-\ell(E)\) divides two for a curve centre.
\end{lemma}

\begin{proof}
Let \(S\to E\) be a resolution. The image of \(N^1(Y)_\R\) in
\(N^1(S)_\R\) contains the restriction of an ample class on \(Y\), whose
square is positive. The Hodge index theorem makes the intersection
pairing nondegenerate on this image. Since every curve on \(E\) is
dominated by a curve on \(S\), the kernel of the restriction map is
exactly the annihilator of the curves on \(E\). This proves the rank
assertion.

The sum \(\operatorname{rank}B_E+E^3\) is monotone along
\(E\)-negative elementary flops. Wilson's formula changes the form
by a term of rank at most one while increasing the integral cube by
at least one \cite[Proposition~2.1]{Wilson}. We can therefore bound
the sum on a model where \(E\) has an elementary contraction. Indeed,
run the program for \((Y,\epsilon E)\), with \(0<\epsilon\ll1\).
Pseudoeffectivity excludes a Mori fibre space. At a nef endpoint,
log abundance would make a multiple of \(E\) semiample, contrary
to rigidity. The preceding small steps are smooth \(E\)-negative
flops \cite{Kollar,Fujino}. The first divisorial step contracts \(E\),
since curves sweeping any other divisor have nonnegative \(E\)-degree.

On the elementary model, pullbacks from the target restrict through
the centre of \(E\), and the relative Picard number is one. A curve
centre therefore gives rank at most two and \(E^3\le8\) by
Lemma~\ref{lem:curve-cube}. A point centre gives rank one and
\(1\le E^3\le9\) by the degree bound for exceptional Gorenstein del
Pezzo surfaces, including the nonnormal case
\cite[\S3, Theorem~5.2 and Remark~5.3]{Gross}.

Write \(d=-\ell(E)\). For a curve centre, \(d\) divides the
anticanonical degree two of the general fibre, as computed in
Lemma~\ref{lem:curve-cube}. For a point centre \(E\) is rationally chain
connected \cite[Corollary~1.4]{HaconMcKernan}. Its resolution is uniruled
and has \(p_g=0\), so its numerical divisor lattice is unimodular by
Lefschetz \((1,1)\) and Poincar\'e duality. The pullback of \(-E|_E\)
pairs with this lattice in \(d\mathbb Z\) and is therefore divisible
by \(d\). Thus \(d^2\mid E^3\), which gives \(d\le3\).
\end{proof}

\begin{proof}[Model independence and the lower bound in Theorem~\ref{thm:fibre}]
The annihilator of \(V_p\) consists of divisors of degree zero on every
curve over \(p\). Consider an elementary flop over \(Z\). If its locus
misses the fibre, the fibre itself is unchanged. Otherwise a rational
divisor in the annihilator has degree zero on the flopping ray, so a
multiple descends to the common small contraction by relative
basepoint freeness. Its pullback on the other side has degree zero on
the flopped curves and on every strict transform of a curve over \(p\).
Reversing the flop identifies the rational annihilators. The spaces
\(V_p\) are rational and any two models are connected by such flops
\cite{KawamataFlops}.

For each exceptional prime whose centre contains \(p\), choose a
sweeping class on a model making that prime elementary over \(Z\), as
in Section~\ref{sec:weights}. If its centre is \(p\), the class belongs
to \(V_p\) by model independence. For a curve centre, resolve the
prime and take the Stein factorization over the centre. Fibres of the
resulting morphism to a smooth curve are numerically equivalent. A
general fibre gives a positive multiple of the sweeping class, while
a fibre over a point lying above \(p\) gives a cycle supported over
\(p\). Thus all \(k+c_p\) chosen classes lie in \(V_p\).

Their intersection matrix with the corresponding primes has negative
diagonal and nonnegative off-diagonal entries. By
Lemma~\ref{lem:row-witness} there is a positive vector that the
matrix sends to a strictly negative vector, since all the classes
annihilate \(D\). Scaling the columns by this vector makes the matrix
strictly diagonally dominant. The chosen classes are independent.
\end{proof}

\subsection{A cone of exceptional divisors}

Write \(E_1,\ldots,E_k\) for the exceptional surfaces over \(p\).
We want an effective combination whose negative is nef on the fibre.
Allowing flops over \(Z\) leads to the following cone of coefficients. Put
\[
 B(x)=\sum_i x_iE_i,\qquad
 \mathcal K=\{x\in\R^k:D-\varepsilon B(x)\in M
                 \text{ for all sufficiently small }\varepsilon>0\},
\]
where \(M\) is the closed movable cone of \(Y\). Near the big class
\(D\), the movable cone and its chamber decomposition are locally
rational polyhedral \cite[Theorems~2.5--2.6 and Corollary~2.7]{Kawamata}.
Thus \(\mathcal K\) is polyhedral with a finite subdivision. On each
piece there is a smooth model over \(Z\) on which \(L=mD-B(x)\) is
nef for \(m\gg0\). Call the model \emph{adapted} to \(x\). The ample
model of a big real nef class is the one defined by a rational point
of the same open nef face.

The outer walls of \(\mathcal K\) come from divisorial contractions,
since at a small wall we can flop to the next chamber. Every
exceptional prime admits a sweeping class over its ample model by the
relative program in Section~\ref{sec:weights}. Such a class has
negative degree on its prime, which in particular is rigid in every
multiple. A general sweeping curve over any other point, or over a
curve centre on \(Z\), misses the \(E_i\). Thus only the classes
\(\ell_a\) sweeping some \(E_i\) impose inequalities. In coefficients
they read
\begin{equation}\label{eq:fibre-cone}
  R_a(x)=-\ell_a(B(x))=d_a x_i-\sum_{j\ne i}a_{aj}x_j\ge0,
  \qquad d_a\in\{1,2,3\},\quad a_{aj}\in\mathbb Z_{\ge0}.
\end{equation}
The diagonal bound is Lemma~\ref{lem:degree}. Applying
Lemma~\ref{lem:row-witness} to all these rows supplies a single
\(v>0\) for which every inequality is strict.

\begin{lemma}[The exceptional fibre and its possible facet diagrams]\label{lem:fibre-cone-finite}
Every nonzero \(x\in\mathcal K\) has strictly positive coordinates.
On every model adapted to this \(x\),
\begin{equation}\label{eq:fibre-cover}
                     f^{-1}(p)_{\rm red}=\bigcup_{i=1}^k E_i.
\end{equation}
The section \(\sum_i x_i=1\) is compact. For \(k\ge2\), retain one row
per facet and draw \(i\to j\) if a retained row for \(E_i\) has
\(a_{aj}>0\). This graph is strongly connected, and for fixed \(k\)
only finitely many labelled facet presentations occur.
\end{lemma}

\begin{proof}
Compare \(x\in\mathcal K\) with the interior vector \(v\) at a
coordinate minimizing \(x_i/v_i=\mu\). For any row for \(E_i\),
\(0\le R_a(x)\le\mu R_a(v)\), so \(\mu\ge0\).
If \(x\ne0\), then \(B(x)\) is effective, meets the fibre, and
\(-B(x)\) is nef on the fibre of an adapted model. A component of the
fibre meeting the support properly would contain a curve of positive
\(B(x)\)-degree, obtained by ample sections on that component.
Consequently each component is contained in the support or disjoint
from it. Connectedness forces the support to contain the whole fibre.
In particular every surface \(E_i\) occurs in \(B(x)\). This proves
strict positivity, \eqref{eq:fibre-cover}, and compactness of the
closed section inside the standard simplex.

For strong connectivity, suppose a nonempty proper set of vertices has
no outgoing arrow. Set their coordinates in \(v\) to zero. The rows
with indices in this set now vanish, and the other inequalities can
only improve. The resulting nonzero point of \(\mathcal K\) contradicts
strict positivity. To see finiteness, an arrow \(i\to j\) gives
\(v_j\le3v_i\). A return path has length at most \(k-1\), so
\(v_i\le3^{k-1}v_j\). The strict facet inequality then gives
\(a_{aj}v_j<d_av_i\) and hence \(a_{aj}<3^k\). The rows are integral,
so there are only finitely many possible labelled presentations.
\end{proof}

If \(k=1\), take \(x_1=1\). The class \(L\) lies off all walls for
\(E_1\), so its ample model preserves \(E_1\) as a surface. Consequently
\(L|_{E_1}\equiv-E_1|_{E_1}\) is big and nef, and \(E_1^3\ge1\).
Equation~\eqref{eq:fibre-cover} and Lemma~\ref{lem:degree} give
\(\dim V_p\le9\). We assume \(k\ge2\) below.

\subsection{Bounding the surface cubes on one model}

For \(k\ge2\) the outer cone belongs to a finite list, but its
subdivision into nef chambers can still vary. To compare the cubics
of all these chambers, consider the loss of volume from \(mD\), given by
\[
 \Phi(x)=\frac{m^3D^3-\operatorname{vol}(mD-B(x))}{6}.
\]
On an adapted model \(L=mD-B(x)\) is nef and \(D|_{E_i}\equiv0\),
so expanding the cube gives \(\Phi(x)=B(x)^3/6\). In particular
\(\Phi\) is independent of \(m\gg0\), homogeneous and piecewise cubic.
Write \(T_i=\partial_i\Phi\), \(H=\nabla^2\Phi\) and \(h_i=H_{ii}\).
On each chamber these derivatives are intersection numbers, namely
\[
   2T_i=L^2E_i,\qquad H_{ij}=B(x)E_iE_j=-LE_iE_j,\qquad
   \partial_i h_i=E_i^3\quad\text{on a chamber}.
\]
The sign and the kernel of the Hessian can both be seen on a surface.
For rational \(x\), take a smooth general member \(P\in|nL|\) of a
basepoint-free multiple and put \(A_i=E_i|_P\). Then
\[
                 H_{ij}(x)=-\frac{1}{n}(A_i\cdot A_j)_P .
\]
All the \(A_i\) are orthogonal to \(D|_P\), whose square is
\(nLD^2=nmD^3>0\). Hodge index makes the negative intersection
pairing an inner product on this orthogonal space of numerical
classes. For \(i\ne j\) a general \(P\) contains no component of
\(E_i\cap E_j\), so the traces have no common component. Thus \(H\succeq0\)
and \(H_{ij}\le0\) for \(i\ne j\). These signs and \(T_i=L^2E_i/2\ge0\)
extend chamberwise to real classes by continuity.

In particular, on the ample model of a rational class we have
\begin{equation}\label{eq:fibre-trace}
       h_i(x)=0 \quad\Longleftrightarrow\quad
       E_i\text{ maps to a point}.
\end{equation}
Indeed a general hyperplane misses a point and meets every
positive-dimensional projective image. In the latter case \(A_i\)
is a nonzero effective curve on \(P\), hence numerically nonzero.

Two facts control the Hessian between chambers and at the outer
boundary. First, for each fixed \(u\) the function
\(q_u(x)=u^{\mathsf T}H(x)u\) is nonnegative, homogeneous of degree
one, and concave. Across an elementary flop let
\(\theta(x)=B(x)C=-LC\) for the old flopping class \(C\), positive
on the new side. Wilson's formula changes \(\Phi\) by
\(-\lambda\theta^3/6\) for some \(\lambda\ge0\)
\cite[Proposition~2.1]{Wilson}. The cubics therefore agree to second
order and \(q_u\) changes by \(-\lambda\theta(u)^2\theta(x)\), a
decrease in slope. A segment in a wall or passing through a face of
higher codimension is a limit of generic parallel segments in the
big movable interior. There one applies the same formula to minus the
volume. Local finiteness and agreement through second order give the
claim on the limiting segment, also at the boundary of \(\mathcal K\).

Second, on a facet labelled by \(E_i\) the ample model contracts
\(E_i\), so \(T_i=L^2E_i/2=0\). Tangential differentiation shows
that the \(i\)-th row of \(H\) is normal to that facet, including on
its closure. Inside every chamber we also have the one-sided bound
\(\partial_i h_i=E_i^3\le9\) by Lemma~\ref{lem:degree}. The following
lemma turns this boundary information into a uniform bound throughout
the coefficient cone.

\begin{lemma}[Bounds on the Hessian diagonals]\label{lem:fibre-diagonals}
Each \(h_i\) is bounded on \(\mathcal K\cap\{\sum x_i=1\}\) by a constant
depending only on the labelled facet presentation of \(\mathcal K\).
\end{lemma}

\begin{proof}[Deduction of the upper bound in Theorem~\ref{thm:fibre}]
Fix an interior point \(b\) with \(\sum b_i=1\) and choose \(s>0\)
such that \(b\pm se_i\in\mathcal K\) for all \(i\). Both choices use
only the labelled cone. Take a full-dimensional chamber whose closure
contains \(b\). The linear piece of the concave function \(h_i\) on
this chamber supports \(h_i\) from above on the whole cone. Its slope
along \(e_i\) is \(E_i^3\). Hence
\[
              0\le h_i(b+se_i)\le h_i(b)+sE_i^3 .
\]
This bounds every cube from below on the same model. Since its
surfaces cover the fibre by \eqref{eq:fibre-cover},
Lemma~\ref{lem:degree} gives
\[
     \dim V_p\le\sum_i\operatorname{rank}B_{E_i}
            \le10k-\sum_i E_i^3\le10k+s^{-1}\sum_i h_i(b).
\]
Lemma~\ref{lem:fibre-diagonals} bounds this expression in terms of
the labelled outer cone. Maximizing over the finite list in
Lemma~\ref{lem:fibre-cone-finite} and rounding up gives \(G(k)\).
Proposition~\ref{prop:two-surfaces} gives the stated small values.

\end{proof}

\begin{remark}
The same supporting-function argument bounds every mixed product on
this model. Put \(C_i=(E_iE_jE_l)_{jl}\). Applying the argument to
\(q_u\) at \(b\pm se_i\) gives
\(0\le q_u(b\pm se_i)\le u^{\mathsf T}(H(b)\pm sC_i)u\), or
equivalently, in the order on symmetric matrices,
\begin{equation}\label{eq:fibre-order}
                      -H(b)\preceq sC_i\preceq H(b).
\end{equation}
In particular \(s|E_iE_jE_l|\le\sqrt{h_j(b)h_l(b)}\). Apply
Cauchy--Schwarz to the \((j,l)\)-entry of each of the positive semidefinite
matrices \((H(b)\pm sC_i)/2\), then to the sum of the two bounds.
\end{remark}

\begin{proof}[Proof of Lemma~\ref{lem:fibre-diagonals}]
Write \(S=\mathcal K\cap\{\sum x_i=1\}\). The facet normals link
different diagonal entries of \(H\). A bound for one entry will
therefore propagate to all of \(S\). We obtain the first bound
from a vanishing row, or from a metric determined by the facets.

\smallskip\noindent\emph{A vanishing diagonal.}
Suppose \(h_i(b)=0\) and \(b+te_i\in\mathcal K\) for all small
\(t>0\). Positive semidefiniteness gives \(H_{ij}(b)=0\) for every
\(j\). By finiteness of the subdivision, a full-dimensional piece has
a closure containing an initial segment of this ray, even if the ray
lies on a wall. Write \(h_i(x)=\sum_j g_jx_j\) on that piece.
Symmetry of its cubic tensor gives \(g_j=\partial_j H_{ii}
=\partial_i H_{ij}\). Thus \(H_{ij}(b+te_i)=t g_j\). For \(j\ne i\)
the left side is nonpositive, so \(g_j\le0\), while
\(g_i=E_i^3\le9\). The supporting inequality and \(x_j\ge0\) give
\begin{equation}\label{eq:incoming-diagonal}
                       h_i(x)\le\sum_jg_jx_j\le9x_i
                               \qquad(x\in\mathcal K).
\end{equation}

\smallskip\noindent\emph{Propagation along the facet graph.}
Suppose \(h_i\) is bounded on \(S\) and \(j\to i\). Choose a row \(a\)
giving this arrow and a normalized relative-interior point \(b\) of
its facet. There the \(j\)-th Hessian row is proportional to the normal,
so \(H_{ji}=-(a_{ai}/d_a)h_j\). Its two-by-two principal minor gives
\[
       (a_{ai}/d_a)^2 h_j(b)^2\le h_j(b)h_i(b).
\]
This bounds \(h_j(b)\). Moving a fixed small distance in direction
\(e_j\) gives an interior point \(z=b+te_j\), with
\(h_j(z)\le h_j(b)+9t\) by the slope bound. Since \(S\) is compact and \(z\) lies in the interior, choose \(c>0\) with
\(z-cS\subset\mathcal K\). Since a nonnegative concave homogeneous
function is superadditive, \(h_j(z)\ge h_j(cx)\) for every \(x\in S\).
This bounds \(h_j\) on \(S\), with choices depending only on the
facets. Strong connectivity propagates the bound to every vertex.

\smallskip\noindent\emph{Compatibility of the facet normals.}
By \eqref{eq:fibre-trace}, a rational point \(b\) on a closed
\(E_i\)-facet gives the vanishing needed in \eqref{eq:incoming-diagonal}
whenever \(E_i\) maps to a point and \(b+te_i\) is feasible for small
\(t>0\). Suppose no such pair \((b,i)\) exists.

Normalize each facet normal \(n_a=R_a/d_a\) to take value one on
its labelled direction \(e_i\). At a rational point \(b\) of that
facet, including its boundary, tangential vanishing gives
\[
                      H(b)(e_i,-)=h_i(b)n_a .
\]
Thus whenever the trace of \(E_i\) on the surface section survives,
its intersections recover the entire normalized normal. In particular
two different \(E_i\)-facets cannot meet where \(+e_i\) is feasible.
If there were two such facets, join rational relative-interior points
and lower the \(i\)-th coordinate along the segment to the graph
\[
                  x_i=\max_{a\text{ for }E_i}
                         \sum_{j\ne i}(a_{aj}/d_a)x_j .
\]
All other labels' inequalities stay strict, since lowering \(x_i\)
relaxes them. At a rational break where the maximizing normal changes,
both facets are active and \(+e_i\) is feasible, a contradiction.

There is therefore one facet per prime (strong connectivity ensures
every label occurs). Their normals form an invertible matrix by strict
diagonal dominance after column scaling by \(v\). The cone is thus
linearly isomorphic to an orthant, so any two facets meet in a ridge.
If \(i\to j\) but \(j\not\to i\), choose a rational relative-interior
point \(b\) of that ridge. Here \(+e_i\) is tangent to facet \(j\)
and points strictly inside the half-space of facet \(i\), so \(h_i(b)>0\). Symmetry of intersection
on the surface section now gives the contradiction
\[
       0> h_i(b)n_i(e_j)=H_{ij}(b)=h_j(b)n_j(e_i)=0 .
\]
Consequently every arrow has its reverse. In both ridge arguments the
obstruction is the same: the prescribed intersections are incompatible
with a nonzero trace of \(E_i\).

\smallskip\noindent\emph{A metric for the surviving traces.}
At a general rational point of its facet \(F_i\), the ample model
contracts \(E_i\) to a curve by our assumption and \eqref{eq:fibre-trace}.
We construct an inner product for which \(e_i\) is normal to \(F_i\);
then \(\partial_i\Phi=0\) is an ordinary normal derivative condition.

Make \(E_i\) elementary over the ample model of a general rational
point of \(F_i\). Its elementary centre must still be a curve, so the
corresponding row \(R_i\) has \(d_i\mid2\) by Lemma~\ref{lem:degree}
and defines the same facet. Thus the normals \(\alpha_i=2R_i/d_i\)
have diagonal two and negative integral entries on the edges. Their
undirected graph is a tree. Indeed, at a vertex minimizing \(v_i\) on
a cycle, the two neighbours \(j,l\) on that cycle would give
\(2v_i>v_j+v_l\ge2v_i\), since \(\alpha_i(v)>0\). Along the tree
choose \(s_i>0\) successively so that \(G_{ij}=s_i\alpha_i(e_j)\) is
symmetric. This symmetrized Cartan matrix is positive definite, since
\(\operatorname{diag}(v)G\operatorname{diag}(v)\) is strictly
diagonally dominant with positive diagonal. By construction,
\(G(e_i,-)=s_i\alpha_i\), so \(e_i\) is normal to \(F_i\).

Take the metric trace \(h=\operatorname{tr}(G^{-1}H)\). It is the sum
of \(q_u\) over a \(G\)-orthonormal basis, so it is concave,
nonnegative and homogeneous of degree one. On a cubic piece adjacent
to an open part of \(F_i\), write \(C\) for the third derivative of
\(\Phi\). The boundary condition \(\partial_i\Phi=0\) gives
\(C(e_i,t,t)=0\) for every tangent vector \(t\) to the facet.
When we take the trace of \(C(e_i,-,-)\), only the normal direction
\(e_i/\sqrt{G_{ii}}\) contributes. Thus on that piece
\[
                 \partial_i h=\frac{C(e_i,e_i,e_i)}{G_{ii}}
                            =\frac{E_i^3}{G_{ii}}\le\frac{9}{G_{ii}} .
\]
Project an interior point \(x\) to \(F_i\) along \(-e_i\). Every other
inequality stays strict. For general \(x\) this projection misses the
lower-dimensional faces of the subdivision of \(F_i\), so the returning
segment starts in a piece with the displayed slope bound. Concavity
can only decrease the slope along the segment, giving the same bound
on \(\partial_i h(x)\). Since
\(x_i>0\), Euler's identity gives
\[
               h(x)=\sum_i x_i\partial_i h(x)
                       \le 9\sum_i x_i/G_{ii}.
\]
Continuity extends the inequality to the closed cone. Since
\(h_i\le G_{ii}h\) by positive semidefiniteness, this bounds all the
diagonals on the section \(\sum_i x_i=1\).
\end{proof}

\subsection{Explicit estimates for small fibres}

For at most three surfaces a direct intersection calculation gives
explicit bounds. We make a surface with smallest coefficient elementary.
On that model all mixed products of three distinct surfaces have a
known sign. With four surfaces, the same argument leaves one product
to check.

\begin{lemma}[An elementary surface of smallest coefficient]\label{lem:smallest-surface}
If \(k\ge2\), one can choose coefficients \(x\in\mathcal K\) and an
adapted smooth model such that the surface with smallest coefficient
is elementary. Its intersection with any two other distinct surfaces
over \(p\) is nonnegative.
\end{lemma}
\begin{proof}
Choose a rational interior vector with distinct coordinates and let
\(x_0\) denote its smallest coordinate, relabelling the surfaces as
\(E_0,\ldots,E_{k-1}\). Decrease only \(x_0\). This relaxes all
inequalities for the other primes in \eqref{eq:fibre-cone}. Since every nonzero
vector of \(\mathcal K\) has strictly positive coordinates, we reach a
wall for \(E_0\) before \(x_0\) vanishes. Rescale the resulting
rational vector so that \(x_0=1\). All other coordinates exceed one.

On an adapted model choose \(m\) large enough that \(L-D\) is nef.
The ample model of \(L\) contracts \(E_0\) and maps to \(Z\), since
\(L=D+(L-D)\) and every \(L\)-trivial curve is \(D\)-trivial. Make
\(E_0\) elementary by flops over this ample model as in
Section~\ref{sec:weights}. Both \(D\) and \(L\) remain nef pullbacks,
so the new model is still adapted.

If the elementary contraction sends \(E_0\) to a point, the numerical
restrictions of divisors to \(E_0\) span a line containing an ample
class. Each \(E_j|_{E_0}\) for \(j>0\) is effective, so its class is a
nonnegative multiple of that ample class. Their pairwise intersections are nonnegative.

If the elementary image is a curve, let \(\ell\) denote the primitive
sweeping class and put \(d=-\ell(E_0)\), \(a_j=\ell(E_j)\). Since \(\ell(B)=0\),
\[
                       d=\sum_{j>0}a_jx_j .
\]
Here \(d\mid2\), the \(a_j\) are nonnegative integers, and every
\(x_j>1\). Consequently exactly one \(a_j\) is positive. For any
other \(j\), the divisor \(E_j\) descends numerically through the
elementary contraction. Its restriction to \(E_0\) is therefore
numerically pulled back from the centre curve. Effectivity forces the
degree on that curve to be nonnegative, by intersecting the restriction
with an ample class on \(E_0\). Thus among the restrictions of any two
distinct \(E_j\), \(j>0\), at least one is nef and the other is effective. Their intersection
is nonnegative.
\end{proof}

\begin{proposition}[Explicit bounds for small exceptional fibres]\label{prop:two-surfaces}
If exactly two exceptional primes map to \(p\), then \(\dim V_p\le12\).
If exactly three map to \(p\), then \(\dim V_p\le262\).
If exactly four map to \(p\), use the model of
Lemma~\ref{lem:smallest-surface} and denote its chosen elementary
surface by \(E_0\). Then \(\dim V_p\le6823\) provided
\(E_1E_2E_3\ge0\) on this model. In particular this holds if the
three surfaces other than \(E_0\) have no common curve there.
\end{proposition}
\begin{proof}
Use the model of Lemma~\ref{lem:smallest-surface} and write \(V(E)\) for
the numerical span of curves on a surface \(E\). The surfaces cover
the connected fibre by \eqref{eq:fibre-cover}. Two meeting surfaces
share a curve, whose numerical class is nonzero by ampleness. Order
the surfaces along a spanning tree of their intersection graph.
Each after the first contributes at most \(\dim V(E_i)-1\) new
directions, so
\[
              \dim V_p\le \sum_i\dim V(E_i)-(k-1).
\]

Order the coefficients as \(1=x_0<x_1\le\cdots\le x_{k-1}\) and
put \(U_i=x_iE_i\). One cubic identity applies to all three cases.
\[
 \sum_i U_i^3
   =\sum_i BU_i^2-\sum_{i<j}BU_iU_j
      +3\sum_{i<j<l}U_iU_jU_l .
\]
The first two sums on the right contribute nonnegatively by the
intersection signs for \(B\). Every triple containing \(E_0\) has
nonnegative intersection by Lemma~\ref{lem:smallest-surface}.
For \(k\le3\) this accounts for every triple, and for \(k=4\) the
only further triple is nonnegative by hypothesis. In particular the
hypothesis holds when \(E_1,E_2,E_3\) have no common curve, since
then \(E_3\) has nonnegative degree on the effective cycle \(E_1E_2\).

For \(k=2\), write \(r=x_1\). From \(\ell(B)=0\) for the elementary
sweeping class of \(E_0\), we have \(r=d/\ell(E_1)>1\), where
\(d=-\ell(E_0)\le3\). Thus \(r\ge3/2\) and the cubic inequality gives
\(E_1^3\ge-E_0^3/r^3\). For a point centre this yields
\(E_1^3\ge-8/3\), hence \(E_1^3\ge-2\), and \(\dim V(E_0)=1\).
For a curve centre \(d\mid2\) forces \(r=2\), so \(E_1^3\ge-1\)
and \(\dim V(E_0)\le2\). In both cases Lemma~\ref{lem:degree}
and the common curve give
\(\dim V_p\le\dim V(E_0)+(10-E_1^3)-1\le12\).

For \(k=3,4\) we use \(\sum_i x_i^3 E_i^3\ge0\) to bound the sum of the unweighted cubes. For each \(j\ge1\)
some arrow of the strongly connected facet graph crosses from
\(\{0,\ldots,j-1\}\) to its complement. Its facet inequality gives
\(x_j\le3x_{j-1}\). Since \(x_1>1\) and each cube is at most nine,
division by \(x_1^3\) gives
\[
 \sum_{i=1}^{k-1} E_i^3
 \ge -\frac{E_0^3}{x_1^3}
     -\sum_{i=2}^{k-1}\bigl((x_i/x_1)^3-1\bigr)E_i^3
 > -9\left(1+\sum_{i=2}^{k-1}(27^{i-1}-1)\right).
\]
By integrality the sum is at least \(-242\) for \(k=3\) and
\(-6794\) for \(k=4\). Lemma~\ref{lem:degree} bounds the contribution
of the last \(k-1\) surfaces by \(10(k-1)-\sum_{i>0}E_i^3\).
Adding at most two for \(E_0\) and subtracting \(k-1\) along the
spanning tree gives the asserted bounds.
\end{proof}

\section{From a movable boundary class to a fibration}\label{sec:reduction}

\begin{proposition}[Bigness in the finite-chamber nonfibered case]\label{prop:reduction}
Suppose that a strict Calabi--Yau threefold \(X\) has finitely many marked nef
chambers and that none of its smooth projective flop models admits a
fibration. If \(\rho(X)\ge19\), every nonzero class in
\(M=\overline{\Mov}^{\,1}(X)\) is big.
\end{proposition}

A nonzero effective rational nef class of cube zero already gives the
desired morphism. Indeed, if \(D\) is an effective representative,
threefold log abundance for \((Y,\epsilon D)\), with \(\epsilon>0\)
small, makes \(D\) semiample. Its numerical dimension is one or two.
A real class on the movable boundary presents two difficulties. Its
ray can be irrational, and even a rational ray in the closure of the
effective cone need not a priori contain a divisor. We deal with these
in turn. The intersection cubic controls the shape of the nef boundary.
Exceptional divisors then provide effective generators for its
supporting faces in \(P=\overline{\Eff}^{\,1}(X)\).

\subsection{Birational walls on a fixed model}

Fix a marked smooth projective model \(Y\). Its nef cone is locally
rational polyhedral on the big locus. Each facet meeting this locus
defines an elementary birational contraction \cite{Kawamata}.
A small contraction gives an adjacent marked flop chamber
\cite{Kollar}. Two distinct chambers can share at most one facet,
since a hyperplane separating their interiors contains their entire
intersection. Thus \(Y\) has finitely many small walls.

\begin{lemma}[Finiteness of divisorial walls on a nonfibered model]\label{lem:div-walls}
If \(Y\) admits no fibration, it has only finitely many elementary divisorial
walls.
\end{lemma}

\begin{proof}
Point contractions and exceptional primes with two elementary rulings
occur only finitely often
\cite[Proposition~1.3 and Remark~1.5]{Uehara}.
If infinitely many walls remain, the proof of
\cite[Proposition~1.10]{Uehara} gives an infinite family of distinct
exceptional primes \(E_i\), each contracted to a curve. Their general
ruling fibres \(r_i\) satisfy
\[
              r_i(E_i)=-2,\qquad r_i(E_j)>0\quad(j\ne i).
\]
Every infinite subfamily contains a nef sum of three of these divisors.
Such a sum is effective and semiample, hence big in the absence of a
fibration.

Let \(H\) be ample. In the affine slice \(H^2E=1\), each \(r_i\)
strictly separates its own normalized class from all the others.
In particular no three of the classes are collinear, since the middle
of three points on a line cannot be separated from both endpoints.
The classes of curves on \(E_i\) span a two-dimensional space \(V_i\), since
pushforward to the centre has one-dimensional image and the elementary
contraction has a one-dimensional kernel. Its projective line
\(\mathbb P(V_i)\) contains \([HE_i]\), so hard Lefschetz and the
noncollinearity just proved allow at most two indices on any one line.
For distinct \(i,j\), the lines meet at \([E_iE_j]\), since the intersection
is effective and nonzero since \(r_i(E_j)>0\).
Retain one index for each distinct line. Pairwise intersecting lines
are concurrent or coplanar. Indeed, three without a common point form
a triangle, and every line meeting all three sides lies in its plane.

If the lines pass through \([v]\), then \(E_iE_j\in\R v\) for \(i\ne j\).
Discard a basis of the divisor span chosen among the \(E_i\). For each
remaining \(i\), expansion in that basis gives \(E_i^2\in\R v\) as well.
Thus all products on the span \(S\) of the remaining divisors belong to
\(\R v\). Take a big nef triple sum \(D\in S\). Since \(D^3>0\), we
have \(Dv\ne0\). If \(x\in S\) annihilates \(v\), then
\(D(xy)=x(Dy)=0\) for every \(y\in S\), so \(xy=0\).
But \((x,y)\mapsto Hxy\) is nondegenerate on \(S\) by Hodge index, since \(HD^2>0\). Thus evaluation on \(v\)
is injective on \(S\), contradicting noncollinearity.

If the lines lie in a plane, let \(U\subset N_1(Y)_\R\) be its
three-dimensional vector space and \(W\) the span of the \(E_i\).
Multiplication by \(H\) embeds \(W\) in \(U\), and noncollinearity gives
\(HW=U\). Also \(W\cdot W\subset U\), since every \(E_iE_j\) is a class
of curves on \(E_i\). The form \(Hxy\) is nondegenerate on \(W\) by
the same big nef triple argument. If \(h\) is the orthogonal projection
of \(H\) to \(W\) for this form, then \(H-h\) annihilates \(HW=U\).
Consequently \(hx^2=Hx^2\) on \(W\), so the plane cubic \(F(x)=x^3\)
has a nondegenerate polar.

Now consider its polar \(E_ix^2\). Every \(x\in\ker r_i\) descends
numerically through the elementary contraction, so its restriction to
\(E_i\) comes from the centre curve. Thus \(E_ix^2\) vanishes on
\(\ker r_i\). This polar is nonzero on \(W\), since the polar at \(h\)
is nondegenerate, and it has \(r_i|_W\) as a linear factor.
We have obtained infinitely many plane cubic polars, each with a
factor strictly separating its own parameter from all other parameters.
The plane cubic lemma in Appendix~\ref{app:polars} excludes this configuration.
\end{proof}

\subsection{Rational points on the nef boundary}

Write \(\rho=\rho(X)\ge19\). Miyaoka's theorem gives \(c_2(Y)\ge0\) on
the nef cone \cite{Miyaoka}. If this functional vanished,
a finite \'etale abelian cover of \(Y\) would give
\(\rho(Y)\le h^{1,1}(A)=9\) by pullback, contrary to \(\rho\ge19\)
\cite[Corollary~7]{CGG}. A rational nef class of positive
\(c_2\)-degree is semiample
\cite[\S1, Fact~(iv) and the discussion after Conjecture~1.2]{Uehara}.
Under our no-fibration hypothesis, every rational nonbig nef class
therefore belongs to \(c_2(Y)^\perp\).

\begin{lemma}[Polyhedrality of the nef, movable and pseudoeffective cones]\label{lem:phases-polyhedral}
Every marked nef chamber is rational polyhedral, as are \(M\) and \(P\).
\end{lemma}
\begin{proof}
The finitely many birational wall hyperplanes cover the big part of
the nef boundary. Take a slice \(H^2L=1\) with \(H\) ample integral.
If some boundary point outside these hyperplanes had positive
\(c_2\)-degree, an entire open patch of the boundary would share these
properties. Each class in the patch is nonbig and thus lies on the
intersection cubic \(F(L)=L^3=0\). The boundary of a convex body is
a topological hypersurface, so this patch cannot lie in the
lower-dimensional singular locus of the reduced cubic. Near a smooth
point of one component, invariance of domain shows that the patch is
open on that component.

We show that every such patch contains a rational direction. The
cubic uses all \(\rho\) variables essentially, since a direction of
constancy would annihilate every product \(xy\), contrary to the
nondegeneracy of \(Hxy\). If \(F\) is geometrically irreducible,
its projective hypersurface is therefore nonconical, of dimension
\(\rho-2\ge17\). Rational points are dense in its smooth real locus
by weak approximation \cite[Theorem~A]{Jones}. If \(F\) is reducible
over \(\C\), it cannot be a product of three linear forms, which would
use at most three variables. Its unique linear factor and remaining
quadratic factor are thus defined over \(\mathbb Q\), up to scaling.
Write \(F=\ell Q\). Since \(\ker\ell\cap\operatorname{rad}Q=0\),
the rank of \(Q\) is at least \(\rho-1\). On the hyperplane \(\ell=0\)
rational density is immediate. On \(Q=0\), the presence of a smooth
real point makes the nondegenerate quotient indefinite. Meyer's
theorem gives a smooth rational point, and projection from this point
gives rational density. In all cases the patch would contain a rational
nonbig nef direction of positive \(c_2\)-degree, a contradiction.

The boundary of the nef cone thus lies in the union of its birational
wall hyperplanes and \(c_2(Y)^\perp\). Their rational supporting
halfspaces cut out the cone, since every segment from an interior point
to an exterior point crosses one of these hyperplanes as it leaves. The finite
union of marked chambers is \(M\), by Kawamata's decomposition on the
big movable locus and passage to its closure \cite{Kawamata}. Since
this union is convex, the chamber generators generate \(M\).

Finally, a prime rigid in every multiple becomes elementary exceptional
after flops, by the negative program in Lemma~\ref{lem:degree}. An ample
rational class in the chamber interior recovers the marked model as
its ample model. Lemma~\ref{lem:div-walls} and the finite number of
chambers therefore give only finitely many rigid primes. For any other
prime \(E\), some \(|mE|\) moves. Its fixed part is supported on \(E\),
since \(mE\) itself is a member. Subtracting that part leaves a
positive movable multiple of \(E\). Hence \(M\) and the finitely many
rigid prime rays generate \(P\).
\end{proof}

\subsection{Effective divisors on the boundary}

The last difficulty is effectivity on the extreme rays of \(P\). A
movable extreme ray has large codimension in \(M\). The order of
vanishing of the volume forces it onto a divisorial wall. The
exceptional prime on this wall will belong to the ray itself.

\begin{lemma}[Numerical dimension bounds the number of nonbig facets through a class]\label{lem:two-outer-facets}
Let \(b\in M\) be nonzero and let \(\nu(b)\) be its numerical dimension
on a marked model where it is nef. At most \(3-\nu(b)\) facets of
\(M\) disjoint from the big cone contain \(b\). If the smallest face of
\(M\) containing \(b\) has larger codimension, then \(b\) lies on a
divisorial wall of a chamber.
\end{lemma}
\begin{proof}
On a chamber through \(b\) expand the volume in the displacement \(x\) as
\[
             (b+x)^3=b^3+3b^2x+3bx^2+x^3.
\]
The first nonzero homogeneous term has degree \(3-\nu(b)\le2\),
since \(bA^2>0\) for any ample \(A\). This term is the same on every
incident chamber. Indeed, chambers through \(b\) are connected by
facet crossings through \(b\), as one sees by taking a generic path
inside \(M\) in a neighbourhood avoiding all nonincident chambers. At each crossing
the cubic changes by a multiple of the cube of the wall equation
\cite[Proposition~2.1]{Wilson}. The change has degree three in \(x\).

Every nonbig facet through \(b\) contains an open piece of a chamber
facet whose closure contains \(b\). The chamber cubic vanishes on
this piece and hence on its supporting hyperplane. Its equation,
which is linear in \(x\), thus divides the common first nonzero term.
There are at most \(3-\nu(b)\) distinct such factors. A face of
codimension \(c\) lies in at least \(c\) facets. If \(c>3-\nu(b)\),
one of its incident facets therefore meets the big cone. A chamber
facet subdividing it and containing \(b\) defines a birational
contraction. A small contraction would give a chamber across this
boundary facet of \(M\), so the contraction is divisorial.
\end{proof}

\begin{lemma}[Supporting faces of a nef class trivial on a divisor]\label{lem:nef-face-restriction}
Let \(b\) be a nonzero integral nef class and \(E\) a prime divisor
on a smooth projective strict Calabi--Yau threefold \(Y\). If
\(b|_E\equiv0\), every linear functional nonnegative on \(P\) and
vanishing on \(b\) also vanishes on \(E\).
\end{lemma}
\begin{proof}
Fix a very ample integral divisor \(A\). The sections of \(A+mb\)
grow without bound by Kodaira vanishing and Riemann--Roch, since the cube
has linear term \(3mA^2b>0\), and all other terms involving \(m\)
are nonnegative. On the other hand, vanishing to any fixed order \(N\)
along \(E\) imposes a bounded number of conditions, uniformly in \(m\).
Indeed the layers of \(NE\) have restrictions \((A+mb-jE)|_E\) for
\(0\le j<N\), each of fixed numerical class as \(m\) varies.
For each layer, choose a very ample integral curve on \(E\) whose
subtraction makes the ample degree negative. Sections then inject
into sections on the normalization of that curve, where the degree
is fixed.

Thus \(A+mb-NE\) is effective for all sufficiently large \(m\)
depending on \(N\). Applying the given functional \(\eta\) yields
\(0\le N\eta(E)\le\eta(A)\) for every \(N\), as required.
\end{proof}

\begin{proof}[Proof of Proposition~\ref{prop:reduction}]
Suppose an extreme ray \(R\) of \(P\) contains no prime. By the
description of \(P\) in Lemma~\ref{lem:phases-polyhedral}, it is a
movable rational ray, also extreme in \(M\). Choose its integral
generator \(b\) and a supporting functional \(\eta\) of \(P\) whose
zero face is \(R\). Then \(\eta\) is strictly positive on every prime
on every marked model.

On a chamber through \(b\), the class \(b\) is nef, nonbig and of
\(c_2\)-degree zero, since positive degree would give a fibration.
As \(R\) has codimension \(\rho-1>2\) in \(M\),
Lemma~\ref{lem:two-outer-facets} places it on a divisorial wall.
Choose a model \(Y\) with this wall, let \(f:Y\to Z\) be its contraction
and \(E\) its exceptional prime. Write \(b=f^*b_Z\) numerically.
If \(f(E)\) is a point, then \(b|_E\equiv0\). If it is a curve \(C\),
the transverse surface estimate of Lemma~\ref{lem:chern}, extended
from ample to nef classes by continuity, gives
\[
               0=c_2(Y)b\ \ge\ \tfrac32\,(b_Z\cdot C)\ \ge\ 0.
\]
Again \(b|_E\equiv0\), since it comes from \(C\). In both cases
Lemma~\ref{lem:nef-face-restriction} puts \(E\) in the zero face of
\(\eta\), a contradiction. Thus every extreme ray of \(P\) contains
a prime.

Finally suppose \(0\ne L\in M\cap\partial P\). On a marked chamber
containing \(L\), express it as a positive real combination of integral
nef ray generators. Each generator that occurs is nonbig, since adding a
pseudoeffective class to a big class preserves bigness. Choose one
such generator. Rational polyhedrality expresses it as a nonnegative
rational combination of the prime generators of \(P\).
Since \(H^1(Y,\mathcal O_Y)=0\), numerical equivalence agrees with
linear equivalence up to torsion for integral divisors. A positive
multiple of this class is therefore effective. Log abundance gives a
morphism with connected fibres. Its image has dimension one or two,
contrary to the hypothesis.
\end{proof}

\begin{proof}[Proof of Theorem~\ref{thm:main}]
For \(\rho(X)\ge19\), Proposition~\ref{prop:reduction} gives
\(M\setminus\{0\}\subset\operatorname{int}P\), so Theorem~\ref{thm:mov}
bounds \(\rho(X)\). Taking the maximum of that bound and \(18\) covers
all ranks.
\end{proof}

\appendix
\section{The cubic of a divisor contracted to a curve}\label{app:curve-cube}

\begin{lemma}[The cube of a divisor contracted to a curve]\label{lem:curve-cube}
If an elementary crepant contraction \(f:Y\to Z\) from a smooth
projective Calabi--Yau threefold maps a prime \(E\) to an integral
curve \(C\), then \(E^3\le8\).
\end{lemma}

\begin{proof}
Adjunction gives \(E^3=K_E^2\). In a general transverse Du Val
resolution the components of the exceptional fibre form a tree.
Monodromy acts transitively on its vertices because \(E\) is
irreducible. Every automorphism preserves the centre of the tree,
which is a vertex or an edge. Thus the general fibre consists of one
rational curve of \((-E)\)-degree two, or two rational curves of
degree one meeting transversely. We work over \(C_0=C^\nu\), of genus
\(g\). Local elementary crepant contractions can have singular curve centres
\cite[Proposition~9.14 and Remark~9.17]{BrownWemyss}.

\smallskip\noindent\emph{An irreducible general fibre.}
Let \(u:S\to E\) be the minimal resolution of the normalization.
Write \(u^*K_E=K_S+B\). The conductor on the normalization is effective
and vertical, since \(E\) is generically a smooth ruled surface over
\(C_0\). If \(B_0\) is its strict transform, then \(B-B_0\) is exceptional.
The divisor \(K_S+B_0\) has nonnegative degree on the resolution curves
by minimality and effectivity. The negativity lemma applied to
\(-(B-B_0)\) gives \(B\ge B_0\). Thus \(B\) is effective and vertical.

Contracting vertical \((-1)\)-curves gives a ruled surface
\(S_0\to C_0\), on which the pushforward of \(B\) has class \(vF\)
for \(v\ge0\) and a fibre \(F\). Since the exceptional intersection
form of \(S\to S_0\) is negative definite,
\[
       K_E^2=(K_S+B)^2\le(K_{S_0}+vF)^2=8(1-g)-4v\le8 .
\]

\smallskip\noindent\emph{A general fibre with two components.}
Let \(\pi:B\to C_0\) be the connected double cover parametrizing the
two lines, and put \(h=g(B)\). The normalization of \(E\) is a ruled
surface \(q:S\to B\). To see this, extend the family of lines over
\(B\) by properness of the relative Hilbert scheme. Every integral
fibre curve \(A\) has \(H^1(\mathcal O_A)=0\), by its ideal sequence,
rational singularities of \(Z\), and \(R^2f_*\mathcal I_A=0\)
(the fibres have dimension at most one). Thus \(A\) is a smooth
rational curve. A flat limit of a degree-one line has a unique curve
component of multiplicity one, since \(-E\) is relatively ample. This
component already has Euler characteristic one, leaving no embedded
or isolated points. The resulting family is ruled, and its evaluation
\(u:S\to E\) is birational. It is finite because it embeds each line and
\(B\to C\) is finite.

The general node defines a section on each line. Properness extends it
to a section \(J\subset S\) and a map \(\gamma:C_0\to E\), with
\(u|_J=\gamma\pi\). First identify conjugate points of \(J\), forming
the pushout \cite[Theorem~5.4]{Ferrand}
\[
        S\xrightarrow{\eta} T=S\amalg_J C_0\xrightarrow{\nu} E .
\]
Functions on \(T\) are functions on \(S\) whose restrictions to \(J\)
descend to \(C_0\). They form a finite intermediate
\(\mathcal O_E\)-algebra. Hence \(\nu\) is finite and birational and is
an isomorphism on the generic fibre. Write \(\sigma\simeq C_0\) for
the image of \(J\).

The canonical class of \(T\) and its differentials along \(\sigma\) can
both be read from the local model.
Over a completed local ring \(R\) of \(C_0\), write the double cover
as \(R'=R[t]/(t^2-a)\). A fibre parameter \(z\) cutting out \(J\) gives
\[
     \widehat{\mathcal O}_T=R+zR'[[z]]
       \simeq R[[x,y]]/(y^2-ax^2),\qquad (x,y)=(z,tz).
\]
Thus \(T\) is Gorenstein, and its section ideal is the conductor of
\(\eta\). On \(S\) this ideal cuts out \(J\). On \(\sigma\) the
corresponding conormal bundle is \(\pi_*\mathcal O_J(-J)\).
The projection to \(C_0\) splits off the base differentials, so
\[
        M:=\Omega_T|_\sigma
          \simeq \Omega_{C_0}\oplus\pi_*\mathcal O_J(-J).
\]

For a finite birational map \(p:W\to W'\) of Gorenstein surfaces,
its conductor on \(W\) is \(\omega_W\otimes p^*\omega_{W'}^{-1}\).
Indeed, finite duality identifies this line bundle with the module
of multipliers \(\operatorname{Hom}_{W'}(p_*\mathcal O_W,\mathcal O_{W'})\)
\cite[Tags 0AX1 and 0FKW]{Stacks}. Write the conductor of \(\nu\) on
\(T\) as \(\mathcal O_T(-D)\). Since \(\nu\) is generically an
isomorphism over \(C_0\), \(D\) is effective and vertical. Set
\(d=\deg D|_\sigma\). Its pullback \(V=\eta^*D\) is a sum of fibres
of \(q\) of total multiplicity \(VJ=2d\).

Applying the conductor formula to \(\eta\) and \(\nu\) now gives
\(u^*K_E=K_S+J+V\). As \((K_S+J)F=-1\) for a fibre of \(q\),
adjunction on \(J\) yields
\[
                 K_E^2=4(1-h)-J^2-4d .
\]
Riemann--Roch for \(\pi\) gives \(\deg M=4g-3-h-J^2\). Eliminating
\(J^2\) gives
\[
                 K_E^2=7-3h-4g+\deg M-4d .
\]

At a general node the local equation of \(E\subset Y\) has zero
differential, so \(\gamma^*\Omega_Y\to M\) is generically an
isomorphism. Both are rank-three bundles. The conductor of \(\nu\)
annihilates the cokernel. Indeed, if \(a\) lies in the conductor and \(b\)
is a function on \(T\), then \(a,ab\) come from \(E\), and
\(a\,db=d(ab)-b\,da\) is an ambient differential with coefficients
on \(T\). Hence
\(M(-D|_\sigma)\subset\gamma^*\Omega_Y\subset M\).
Taking determinants and using \(K_Y\simeq\mathcal O_Y\) gives
\(\deg M\le3d\). We obtain \(E^3\le7-3h-4g-d\le7\). For a smooth centre see also \cite[Proposition~2.4]{Uehara}.

If \(C\) is singular, this improves to six. Indeed,
\(f_*\mathcal O_Y(-E)=\mathcal I_C\) and relative vanishing gives
\(R^1 f_*\mathcal O_Y(-E)=0\). Thus
\((f|_E)_*\mathcal O_E=\mathcal O_C\). If \(D=0\), then
\(E=T\) and the map to \(C\) would factor through its normalization,
contradicting this equality. Hence \(V\) contains a fibre, so
\(d=VJ/2\) is a positive integer.
\end{proof}

\section{Polynomial bounds in the face argument}\label{app:cube-rate}

\subsection{Recognizing a cube from the faces at one vertex}

We prove the quantitative assertion of Lemma~\ref{lem:large-cube}.
Keep \(q=10\) fixed. The purpose of the count is to find \(t\) edge
directions at one vertex such that every \(q\) of them span a cube.
The incidence criterion will then assemble these small cubes into
a single \(t\)-cube.

At a vertex of the simple \(n\)-polytope \(Q\), each choice of \(q\)
directions spans a unique \(q\)-face. We claim that some vertex has
at most an \(O_q(1/n)\) fraction of choices spanning a noncube. Indeed,
the count in Lemma~\ref{lem:large-cube} gives \(O_q(Vn^{q-1})\)
noncubical \(q\)-faces, where \(V\) denotes the number of vertices.
Every cubical \(q\)-face gives exactly \(2^q\) vertex incidences. By the
height estimate, the number of incidences with cubes is at least
\[
       2^q\bigl(f_q-O_q(Vn^{q-1})\bigr)
          \ge V\binom nq-O_q(Vn^{q-1}).
\]
Here \(2^q\binom{\lfloor n/2\rfloor}q\) and \(\binom nq\) have the
same leading term. Comparing with the \(V\binom nq\) incidences
in total proves the claim.

At such a vertex \(v\), choose \(2t\) directions uniformly at random.
The expected number of noncubical \(q\)-subsets is
\(O_q(\binom{2t}q/n)<t\) whenever \(n>C_qt^{q-1}\), enlarging
\(C_q\) so that \(n\ge2t\). For some choice we may delete one direction
from each bad subset and retain \(t\) directions with the required
property.

Let \(F\) be the \(t\)-face spanned by these directions, and \(A\) the
\(n-t\) facets containing \(F\). Along the \(i\)-th edge from \(v\),
one facet \(a_i\) through \(v\) is exchanged for a new facet \(b_i\).
On any of our \(q\)-cubes containing that edge these are its opposite
sides in direction \(i\). In particular, the \(b_i\) are distinct,
because each pair of directions occurs in one such cube.

Choose either \(a_i\) or \(b_i\) for each \(i\). Together with \(A\)
this gives \(n\) facets. Any selection of at most \(q\) of them
involves at most \(q\) directions. Extend those directions to \(q\)
of the retained ones. Their cube lies in all facets of \(A\) and has
a vertex on the prescribed sides. By Lemma~\ref{lem:ten-labels} all
\(n\) facets meet, hence meet at a vertex of \(F\).

Simplicity makes the \(2^t\) resulting vertices distinct. Choices that
differ once share \(n-1\) facets and hence give adjacent vertices.
This supplies all \(t\) edges of \(F\) at every listed vertex.
Connectedness of the graph now gives all vertices, with exactly the
facet incidences of a cube.

\subsection{Repeating a part of a minimal nonface}\label{app:repeated-part}

We justify the extraction used in Proposition~\ref{prop:blocks}. On
\(K_0\) and a family of \(N\) disjoint sets \(D_i\), take the minimal
nonfaces of least pair \((s,k)\). Recall that \(s\) is their largest
intersection size with a \(D_i\), and \(k\) counts the intersections
of this size. Each \(K_0\cup D_i\) is a minimal nonface, so a least
pair exists on every nonempty subfamily.

In a nonface of least pair, mark any one of its largest intersections
and call the remaining labels its remainder. If the same remainder
occurs with marked parts in \(m\) distinct sets, it supplies the
\(K,P_1,\ldots,P_m\) used in the proposition. Otherwise consider
nonfaces of least pair meeting exactly \(r\) sets. A remainder meets
\(r-1\) sets, so there are \(O_q(N^{r-1})\) possibilities for it. Once those sets are chosen, the number of ways
to select labels in them and \(K_0\) is bounded in terms of \(q\). Each remainder extends in fewer than \(m\) further
sets, with at most \(2^q\) choices of the marked part per set. This
bounds the number of these nonfaces by \(O_q(mN^{r-1})\).

In a uniformly chosen subfamily of \(T\le N\) sets, a given nonface
meeting \(r\) sets survives with probability at most \((T/N)^r\).
Since some set supplies at least two labels, \(r\le q-1\). The
expected number of surviving nonfaces of least pair is consequently
at most
\[
                  \frac{C_q m}{N}\sum_{r=1}^{q-1}T^r .
\]
For \(N>C'_q m(T+1)^{q-1}\) this is less than one, so a subfamily
of \(T\) sets avoids that pair. Its least pair is strictly larger.

There are only boundedly many pairs, since
\(2\le s\le q\) and \(1\le k\le\lfloor q/s\rfloor\).
Starting with a target size of one and composing
\(T\mapsto\lceil C'_q m(T+1)^{q-1}\rceil+1\) once for each possible
pair ensures that the family stays nonempty through that many
restrictions. Its least pair cannot increase that many times.
A repeated remainder must therefore occur. Because \(q\) is fixed,
the required initial size is bounded by a polynomial in \(m\).

\clearpage
\section{Further results on Picard rank, fibrations and flops}\label{app:project-results}

By Theorem~\ref{thm:main}, sufficiently large Picard number forces a
fibration on a flop model or infinitely many marked chambers. The results
below concern both possibilities. We use the convention on strict Calabi--Yau
threefolds from Section~1. In particular, the Picard number is unchanged under
flops. Fibered threefolds in this class have bounded Picard number by
Engel--Filipazzi--Greer--Mauri--Svaldi \cite[Corollary~1.2]{CFibered}.

\subsection{Recognizing a fibration from a real nef class}

Let \(Y\) be a smooth marked model and \(D\ne0\) a real nef class on \(Y\)
with \(D^3=0\). Write \(V_D\) for the smallest real linear subspace of
\(N^1(Y)_\R\) defined over \(\mathbb Q\) that contains \(D\). The square
\(D^2\) is a numerical curve class, and thus a linear functional on divisors.

\begin{proposition}[Positive second-Chern degree]\label{prop:c-real-nef}
Suppose \(c_2(Y)D>0\).
If \(D^2=0\), the ray of \(D\) is rational. Its primitive integral generator
\(P\) defines a K3 fibration over \(\mathbb P^1\) by a multiple, and
\(c_2(Y)P\in\{12,24\}\).
If \(D^2\ne0\), one has the equivalences
\[
 D^2|_{V_D}=0\quad\Longleftrightarrow\quad
 \R_{>0}D^2\text{ is rational}\quad\Longleftrightarrow\quad
 D\equiv f^*A
\]
for a genus-one fibration \(f:Y\to S\) and an ample real class \(A\) on \(S\).
\end{proposition}

\begin{proof}[Proof sketch]
Fix a smooth very general very ample surface \(H\). Suppose first that
\(D^2=0\) and its ray is irrational. There is a \(v\in V_D\) with \(HDv>0\).
Otherwise \(D\) lies in the radical of the rational form \(Hxy|_{V_D}\),
which has dimension at most one by the Hodge index theorem.
The line \(tD\) is dense modulo the lattice of \(V_D\), so we can find
integral \(L_j=t_jD+e_j\) with \(t_j\to\infty\), \(e_j\to0\), and
\(e_j/\|e_j\|\to v/\|v\|\).

The restrictions \(L_j|_H\) are eventually nef and big. Indeed their
squares are positive. A curve of bounded \(H\)-degree contradicting
nefness would, on a subsequence, have constant numerical class orthogonal
to \(D\), hence to \(V_D\). For curves of unbounded degree, normalize by
that degree and take a limit on \(H\). Adjunction makes its square
nonnegative, and its pairing with \(D|_H\) vanishes. The surface index
theorem makes it a positive multiple of \(D|_H\), contradicting \(HDv>0\)
and the negative pairing with \(L_j\). Vanishing applied to
\(0\to L_j\to L_j+H\to K_H+L_j|_H\to0\) gives
\(H^2(Y,L_j)=0\), and Serre duality gives \(H^3(Y,L_j)=0\). Hence
\[
 h^0(Y,L_j)\ge t_j c_2(Y)D/12+o(t_j),\qquad H L_j^2=o(t_j).
\]
Restriction of sections to \(H\) is injective since \((L_j-H)DH<0\).
If the system had image of dimension at least two, its restriction to
\(H\) would map onto a nondegenerate surface of degree at least
\(h^0(Y,L_j)-2\) and at most \(HL_j^2\). Hence it is a pencil for large
\(j\). The divisor of a general connected fibre has bounded \(H^2\)-degree since
its coefficient in the moving class is at least \(h^0(Y,L_j)-1\).
Pass to a subsequence of constant integral fibre class \(F\).
Nefness gives \(HDF=0\) from \(HDL_j\to0\), and movability gives
\(HF^2\ge0\). Hodge index forces \(F\) to be proportional to \(D\),
contradicting irrationality.

The integral positive-\(c_2\) semiampleness criterion now applies
\cite[Lemma~2.2(3)]{CLiuSvaldi}. For a finite fundamental group one may
apply it on the simply connected cover and descend by the norm. In the
remaining case the Beauville--Bogomolov cover is abelian or a product
of a K3 surface and an elliptic curve. Positive \(c_2\)-degree excludes
the former, and on the product the square-zero pullback comes from a
positive-degree class on the elliptic factor. The base curve has genus
zero since \(q(Y)=0\). A smooth fibre of class \(mP\) has trivial
canonical bundle and Euler number \(m c_2(Y)P>0\), so it is a K3
surface. Its Euler number is \(24\), and Riemann--Roch gives
\(12\mid c_2(Y)P\).

For \(D^2\ne0\), a pullback as stated has square proportional to the
general fibre. A rational square ray has a rational kernel containing
\(D\), hence containing \(V_D\). Conversely, suppose \(D^2|_{V_D}=0\).
The restriction \(D|_H\) is big and nef and remains so under small
perturbations in \(V_D\), as one sees by writing it as ample plus effective on \(H\).
Only finitely many components of the effective part can obstruct
nefness, and any one orthogonal to \(D\) annihilates \(V_D\).
For integral \(L=tD+e\) with \(e\in V_D\) small and \(t\) large,
the same vanishing argument and Riemann--Roch give
\[
 \frac{\chi(Y,L)}t=\frac{c_2(Y)D}{12}+\frac{De^2}{2}+O(t^{-1})>0.
\]
Serre duality gives \(H^3(Y,L)=0\), so these classes are effective.
Lattice approximation lets us choose errors surrounding zero. A positive
relation among them expresses \(D\) as a positive sum of effective classes.
Choose a small positive multiple of the resulting effective representative
to obtain a klt boundary. The rational polytope of nef boundaries on its support
\cite[6.6]{CFujinoPoly} and threefold log abundance \cite{Fujino} express
\(D\) as a positive sum of semiample rational classes. The Stein
factorization of their product morphism gives \(f\). The base has
dimension two since \(D^2\ne0\), \(D^3=0\), and adjunction gives genus one.
\end{proof}

\subsection{Contractions in high rank}

\begin{proposition}[A target of Cartier rank eighteen]\label{prop:c-contractions}
If \(Y\) has no genus-one fibration and \(\rho(Y)\ge19\), it admits a
birational contraction to a normal projective variety of Cartier Picard
number \(18\). More precisely, starting with any birational contraction
\(f:Y\to T_d\) with \(d=\rho(T_d)\ge19\), one can successively contract
\(T_d\to T_{d-1}\to\cdots\to T_{18}\), with \(\rho(T_i)=i\).
Here the ranks on the targets count Cartier divisors modulo numerical
equivalence.
\end{proposition}

\begin{proof}[Proof sketch]
On \(N^1(T_d)_\R\) write \(F(x)=x^3\) and \(\ell(x)=c_2(Y)f^*x\).
The classification of birational \(c_2\)-contractions gives
\(\ell>0\) on the ample cone for \(d\ge19\)
\cite[Theorems~3.3--3.4]{COguisoSakurai}. If the entire nef boundary
had cube zero, it would contain a relatively open patch of the smooth
locus of \(F=0\) with \(\ell>0\). Indeed \(\ell\) is positive somewhere
on the boundary, and Hodge index at an ample class gives a nondegenerate
polar. Thus the cubic is nonconical and squarefree, and its singular
locus has codimension at least two. At a smooth point the nef boundary
is open in the real cubic by invariance of domain.

For a geometrically irreducible cubic in at least \(19\) variables, weak
approximation gives a rational point of the patch \cite{Jones}.
In the reducible case the nondegenerate polar forces a rational linear
factor times a rational quadratic form of rank at least \(d-1\).
The linear component has dense rational points, as does the quadratic
component at a smooth real point by Meyer's theorem and projection.
Pulling back gives a rational nef class of positive \(c_2\)-degree,
nonzero square and zero cube, hence a genus-one fibration by
Proposition~\ref{prop:c-real-nef}.

There is therefore a big nef boundary point. Local rational
polyhedrality on \(Y\), restricted to \(f^*N^1(T_d)\), gives a facet
meeting the big locus \cite{Kawamata}. A rational class in its relative
interior gives a base-point-free contraction on \(Y\) factoring through
\(T_d\). Rational big classes in that relative interior have the same
contracted curves and give the same target. They span the facet, while
every pullback from the target annihilates its contracted curves. The
target therefore has Cartier rank exactly \(d-1\), and the argument
iterates. The relative Picard number of the composite can still be large.
\end{proof}

\begin{proposition}[Bounded integral polarizations]\label{prop:c-polar}
Fix \(v>0\). For all strict Calabi--Yau threefolds \(Y\) and integral
big and nef \(L\) with \(L^3\le v\), the divisor lattice modulo torsion,
with its integral form \(B_L(D,E)=LDE\), has only finitely many isometry
classes. In particular \(v\) bounds \(\rho(Y)\), including directions
in the radical of \(B_L\).
\end{proposition}

\begin{proof}[Proof sketch]
The system \(|5L|\) gives a base-point-free birational map
\cite[Theorem~2 and the following remark]{CJieLiu}. Its nondegenerate
image has degree \(125L^3\) and hence bounded ambient dimension.
These images form a bounded family. On its normalized image \(T\), let
\(A\) be the hyperplane bundle. The contraction \(f:Y\to T\) satisfies
\(f^*A=5L\), and the volume form shows that \(T\) is Gorenstein canonical
and \(f\) crepant. Stratify the bounded family and spread out a resolution
of each generic fibre, separating exceptional components after finite
base change. Smooth proper local triviality gives finitely many lattices
\(H^2(W,\mathbb Z)/\mathrm{tors}\) on these resolutions \(\pi:W\to T\),
with their \(\pi^*A\)-pairings and marked exceptional classes. On fibres
birational to \(Y\) these are divisor lattices by the exponential sequence.

Every prime divisor on \(Y\) has a strict transform on \(W\), since a
divisor exceptional over the smooth \(W\) would have positive
discrepancy over \(T\). The induced pushforward on divisor classes
gives the lattice of \(Y\) as a quotient by the saturation \(K\) of the
classes of the exceptional primes with positive discrepancy. To recover
the pairing on this quotient, cut \(T\) by a general hyperplane. Its
pullback to \(Y\) is the minimal resolution of the normal surface section
by adjunction, so the section on \(W\) maps to it. On a common resolution \(W\xleftarrow{p} U\xrightarrow{q}Y\), the
integral lift \(p_*q^*D\) restricts to the pullback of the class on the
minimal surface resolution. It is therefore orthogonal to \(K\), whose
generators restrict to zero or to exceptional curves. The pairing of
two such lifts is \(5B_L\), independently of the orthogonal lifts chosen. Thus the finitely many marked resolution lattices and choices
of exceptional primes determine the integral form even when it is
degenerate.
\end{proof}

\subsection{A sharp bound for complete intersections}

\begin{proposition}[Picard rank at most two]\label{prop:c-cicy}
Suppose \(X\) is a scheme-theoretic complete intersection of
multihomogeneous hypersurfaces in a product of projective spaces.
If no smooth projective flop model admits a fibration, then
\(\rho(X)\le2\). Equality occurs.
\end{proposition}

\begin{proof}[Proof sketch]
Every nonconstant projection has three-dimensional image, since a
smaller image gives a fibration by Stein factorization. Thus its
hyperplane class \(H_i\) is big and nef. Constant factors of total
dimension \(q\) can be removed since their conormal directions contribute
a trivial summand of rank \(q\), whereas a negative defining
multidegree involving a nonconstant factor has no sections on \(X\).
Consequently \(q\) independent equations cut out the constant point
locally, and substitution gives a complete intersection in the remaining
product. Let \(r\) be its number of factors. Vanishing for the negatives
of the \(H_i\) and of the defining multidegrees, in the restricted Euler
and conormal sequences, identifies \(H^1(\Omega_X^1)\) with
\(\bigoplus_i\C c_1(H_i)\). In particular \(r=\rho(X)\).

Let \(s_i\) count equations supported on factor \(\mathbb P^{n_i}\)
alone and put \(d_i=n_i-s_i\ge3\). The differentials of these equations
are independent on the projected image. Adjunction gives total degree
\(n_i+1\) on factor \(i\). The \(\sum d_i-3\) mixed equations each use
at least two degrees, leaving
\[
        2\big(\textstyle\sum d_i-3\big)\le\sum(d_i+1),
        \qquad 3r\le\sum d_i\le r+6 .
\]
Thus \(r\le3\). If \(r=3\), every pure equation is linear, and eliminating
them leaves two bilinear equations on each pair in \((\mathbb P^3)^3\).

Write the six equations as \(x^t A_i y=x^t B_i z=y^t C_i z=0\)
for \(i=1,2\). For general coefficients the four columns
\(A_1y,A_2y,B_1z,B_2z\) have one relation at a general point.
Its two pairs of coefficient ratios give a rational map to
\(\mathbb P^1\times\mathbb P^1\). Fix general ratios \((a,b)\).
Writing \(A_a=a_1A_1+a_2A_2\), and similarly for \(B_b\), the relation
determines \(y\sim A_a^{-1}B_bz\). The last two equations now give two
general quadrics in \(\mathbb P^3_z\), since bilinear forms on the graph of an
isomorphism restrict onto all quadrics. On a dense open of their
intersection three columns are independent, so their annihilator
recovers \(x\). This gives a genus-one generic fibre. The specialization
inequality for fibering genus \cite[Proposition~3.6]{CChen} extends the
conclusion to each smooth member, since a smooth \(K\)-trivial threefold
is nonuniruled.

On a resolution of this rational fibration the discrepancy divisor is
effective and has degree zero on a general elliptic fibre. All its
components are therefore vertical. Pushing down a very ample system
from the base gives a movable divisor whose system has image of
dimension at least two. Its degree on the covering fibres is zero, so
it is non-big and its Iitaka dimension is exactly two. Directed flops and log abundance
make it semiample. This excludes \(r=3\). For sharpness, Oguiso's
general complete intersection of degrees \((1,1),(1,1),(2,2)\) in
\(\mathbb P^3\times\mathbb P^3\) has rank two and two irrational
outer movable rays \cite[Proposition~6.1]{COguiso}. A fibration on any
flop model would give a nonzero rational movable class of zero volume,
necessarily on one of these rays.
\end{proof}

\subsection{Specialization of a fibration}

\begin{proposition}[Preserving the base dimension]\label{prop:c-specialization}
Let \(\mathcal X\to T\) be flat and projective over a smooth complex curve,
with normal integral total space and geometrically integral generic
threefold. Suppose a multiplicity-one component of a special fibre is
birational to \(X\). If the geometric generic threefold admits a rational
fibration over a base of dimension \(d\in\{1,2\}\) whose smooth connected
general fibre has Kodaira dimension zero, some smooth flop model of \(X\)
admits a fibration over a base of dimension \(d\). When the generic
fibration is defined over \(\C(T)\), any multiplicity is allowed.
\end{proposition}

\begin{proof}[Proof sketch]
Descend the generic fibration by a finite base change. At a component
of multiplicity one the original parameter is a uniformizer, so the
ramified change of parameter is totally ramified there and preserves
the residue field \(\C(X)\). Without base change this step is unnecessary.
Restrict the valuation of the component to the function field \(L\)
of the total base of the fibration. Its restriction \(w\) is nontrivial.
The two relative Abhyankar inequalities bound
\(\operatorname{trdeg}_{\C}\kappa(w)\) by \(d\) and
\(\operatorname{trdeg}_{\kappa(w)}\C(X)\) by \(3-d\).
Their sum is exactly three, so the induced rational base has dimension \(d\).

Realize \(w\) as a divisor on a smooth base model and resolve the map,
retaining the component representing \(\C(X)\). A general curve transverse
to the base divisor reduces the comparison of plurigenera to a
one-parameter degeneration. Takayama's inequality
\cite[Theorem~1.2]{CTakayama} bounds the plurigenera of the resulting
general fibre components by those of the original fibre, hence by one.
The volume form on \(X\), divided locally by a base volume form, gives
a nonzero canonical section on the general fibre. After Stein
factorization its Kodaira dimension is therefore zero.

Resolve the induced rational map as \(X\xleftarrow{p}W\xrightarrow{q}B\).
For general very ample \(A\) on \(B\), the movable \(D=p_*q^*A\) has
\(p^*D=q^*A+G\), where \(G\) is effective and exceptional.
Smoothness of \(X\) gives \(G\le a K_W\) for some \(a\), using the
canonical divisor of its volume form. On a general fibre \(F\) this bounds
the sections of \(mp^*D\) by those of \(maK_F\), so the complete systems
of \(D\) are constant along that fibre. The subsystem from \(B\) has
image dimension \(d\), whence \(\kappa(D)=d\). A directed minimal model
program consists of flops by movability, and log abundance gives the
asserted morphism.
\end{proof}

\subsection{Flop walls and the second Chern class}

Set \(V=N^1(X)_\R\), \(M=\overline{\Mov}^{\,1}(X)\),
\(C=\operatorname{int}M\) and \(\Omega=\mathbb P(C)\).
Write \(k_Y(u)=c_2(Y)u\), taking the Chern class on the marked model \(Y\).

\begin{proposition}[Continuation and a crossing bound]\label{prop:c-walls}
If a facet between two marked nef chambers spans \(H\subset V\), the
closed flop facets in \(H\) cover \(H\cap C\). These hyperplanes form a
locally finite arrangement in \(C\). Moreover, for a primitive integral
\(v\in M\) and any smooth marked model \(Y\), some marked model on which
\(v\) is nef is reached from \(Y\) after at most \(\lfloor k_Y(v)/2\rfloor\)
elementary flop wall crossings.
\end{proposition}

\begin{proof}[Proof sketch]
At an interior movable class \(u\), choose a model on which it is nef.
Choose a rational big class in the relative interior of its minimal nef
face and take its contraction.
A divisorial contraction is impossible, since its exceptional divisor
would be swept by
contracted curves, whose classes are nonnegative on \(M\) and vanish at
its interior point \(u\). For this small contraction \(Y\to Z\), degrees
on the exceptional curves give an injection
\[
  N^1(Y/Z)_\R\lhook\joinrel\longrightarrow
       \bigoplus_{z} N^1(Y_z/Z_z)_\R
\]
into the divisor spaces of the exceptional germs. Locally the nef
chambers are the chambers of finite restricted ADE root arrangements
\cite{CNabijouWemyss}. Restrict their product to the image of the
injection. For a rational class in any open cell, twisting by a
pullback from \(Z\) gives a movable effective representative because the
contraction is an isomorphism away from its exceptional curves.
A relative directed minimal model program flops it to a nef model.
The chosen class avoids the local walls, so it is ample at every
exceptional fibre, hence relatively ample. Since contracted curves
test relative nefness, the entire cell is precisely the chamber of
this model. This realizes exactly the local chambers meeting the global
divisor space. Kawamata's local chamber description
\cite[Theorem~2.6]{Kawamata} gives the corresponding arrangement near \(u\).
The locus in \(H\cap C\) where \(H\) occurs as a wall is open by this
description and closed by local finiteness. The section is connected.

For the count, orient a primitive integral curve normal \(\eta\) to be
positive on the departing ample chamber. Wilson's formula
\cite[Proposition~2.1]{Wilson} changes \(k_Y\) by
\(2(\sum_{j\ge1} j n_j)\eta\) with \(n_j\) nonnegative integers,
at least one positive. To see that \(k_Y\) is nonnegative on \(M\), direct
a rational big movable class by flops to a nef model, apply Miyaoka's
nonnegativity there, and then pass to the closure. Follow the segment
from an ample point of \(Y\) towards \(v\), choosing the starting point so
the open segment avoids all faces of codimension at least two. At each crossing
\(\eta(v)\le-1\), so the nonnegative integer \(k_Y(v)\) drops by at
least two. After the last crossing the endpoint is in the closed nef
chamber. This also applies to a rational ray of zero volume.
\end{proof}

\begin{proposition}[Two square-zero classes on different models]\label{prop:c-vol}
Let \(D\) and \(E\) be real nef classes on marked models \(Y\) and \(Z\),
respectively, with \(D_Y^2=E_Z^2=0\). Put
\(B=(k_Z-k_Y)(D-E)\) in the common marking. Then
\[
             B\ge0,\qquad \operatorname{vol}(D+E)\le B^3/32 .
\]
\end{proposition}

\begin{proof}[Proof sketch]
Perturb the endpoints to general ample \(d,e\) on their respective
models and follow their segment. At the \(i\)-th crossing put
\(x_i=\eta_i(d)>0\), \(y_i=-\eta_i(e)>0\). The cubic and Chern jumps
are \(-a_i\eta_i^3\) and \(b_i\eta_i\) with
\(a_i=\sum_j n_jj^3\le(b_i/2)^3\), \(b_i=2\sum_j n_jj\).
Thus \(\sum_i b_i(x_i+y_i)=(k_Z-k_Y)(d-e)\).
Expand the cubic at \(d+e\) on the chamber containing its ray, that
of the midpoint. Transport the terms \(e^3+3de^2\) from the ending
model to the starting model, then add the jumps up to the midpoint.
The contribution at a wall is
\[
 a_i\big(3x_i y_i^2-y_i^3+(y_i-x_i)_+^3\big)
 =a_i m_i^2(3M_i-m_i)
 \le \tfrac14 a_i(x_i+y_i)^3 ,
\]
where \(m_i=\min(x_i,y_i)\), \(M_i=\max(x_i,y_i)\) and
\(t_+=\max(t,0)\). The four endpoint terms tend to zero by the
square-zero hypotheses. Sum the estimates and use continuity of volume.
\end{proof}

\subsection{When the birational group acts cocompactly}

Let \(\Gamma\) denote the image of \(\operatorname{Bir}(X)\) on the
integral numerical divisor lattice. Say that \(\Gamma\) \emph{divides}
\(\Omega\) if it acts properly discontinuously and cocompactly on the
whole projectivized movable interior.

\begin{proposition}[A divided movable cone]\label{prop:c-group}
Suppose \(c_2(X)\) is numerically nonzero and \(\Gamma\) divides
\(\Omega\). The domain \(\Omega\) is strictly convex. The group
\(\operatorname{Bir}(X)\) is word-hyperbolic, acts properly and
cocompactly on a \(\mathrm{CAT}(0)\) cube complex, and is virtually
compact special. Its virtual cohomological dimension is \(\rho(X)-1\).
If \(C\) is a symmetric cone, then \(\rho(X)\le4\).
\end{proposition}

\begin{proof}[Proof sketch]
In rank one the domain is a point and the polarized automorphism group
is finite. Assume henceforth that \(\rho(X)\ge2\). On a nef chamber put
\(s(u)=k_Y(u)\). The wall formula makes \(s\)
continuous and \(\Gamma\)-invariant. It is positive on \(C\) because every
\(k_Y\) is a nonnegative integral covector on \(M\) and is nonzero.
Indeed, a model with numerically trivial \(c_2\) admits a finite
\'etale abelian cover \cite{COguisoSakurai}, hence has no flops.
On the slice \(\operatorname{vol}=1\), choose a compact set of orbit
representatives. A ball of uniform radius about each representative lies
in \(C\), and nonzero lattice covectors have norm bounded below. Thus
for constants \(a,b>0\) one has
\[
  \xi(u)\ge a\operatorname{vol}(u)^{1/3},\qquad
  s(u)\le b\operatorname{vol}(u)^{1/3}\quad(u\in C)
\]
for every nonzero integral covector \(\xi\) nonnegative on \(M\).
To see that \(C=\BigCone(X)\), observe that the discrete integral projective group acts
properly on the properly convex projectivized big cone, so the image
of \(\Omega/\Gamma\) in that connected Hausdorff quotient is both open
and compact. In particular volume vanishes on \(\partial C\).

Take \(0\ne v\) nef on \(Y\) and \(A\) ample there. The cube of
\(A+tv\) tends to infinity since \(A^2v>0\). The first inequality
with \(\xi=k_Y\) forces \(k_Y(v)>0\). The second, divided by \(t\),
then forces \(v^3>0\). Thus the projectivized closed nef chamber
lies compactly inside \(\Omega\). It is a rational polytope by local
polyhedrality. There are finitely many chamber and wall orbits by
local finiteness and compactness of the quotient.

To see strict convexity, suppose an open triangle with vertices
\([v_i]\) lies in \(\Omega\) with all three edges on its boundary.
Write \(G(t)=\operatorname{vol}(\sum_{i=1}^3t_iv_i)\).
Along a general line the third derivative of volume has only upward
jumps by Wilson's cubic formula, so the fourth derivative is
a nonnegative distribution, also along limiting lines. In a single
coordinate the resulting positive function \(g\) satisfies \(g(0)=0\)
and \(g(t)=o(t^3)\) by homogeneity and continuity at the zero-volume
vertex. Since \(g''\) is convex, a positive slope of
\(g''\) would force cubic growth, and a negative value of the resulting
nonincreasing \(g''\) would force \(g\) eventually negative. Hence
\(g''\ge0\), so \(g(t)/t\) is nondecreasing. It follows that
\(G/(t_1t_2t_3)\) is nondecreasing in every variable and homogeneous
of degree zero, hence constant. One marked chamber meets the triangle
on a relatively open patch. Its cubic consequently agrees with volume
on the whole triangle. Equality with a marked cubic at a big movable
point implies nefness on that model, since outside its nef cone a directed
flop strictly increases the cube, and later steps only increase it.
The whole triangle would lie in one compact chamber, an impossibility.
The triangle criterion and hyperbolicity theorem for divisible domains
now apply \cite[Theorem~4.1 and Corollary~9.1]{CBenoist}.

By Proposition~\ref{prop:c-walls}, the complete wall sections give
halfspaces in \(\Omega\). Their stabilizers act cocompactly on the
sections by local finiteness and compact division. They are therefore
quasiconvex codimension-one subgroups in the hyperbolic Hilbert geometry.
The chamber diameters have a common bound \(L\). If \(m\) walls separate
two general orbit points, their segment has Hilbert length at most
\((m+1)L\). The resulting wall distance is proper, and the finite
collection of wall orbits gives a proper cocompact action on the dual
cube complex \cite{CBergeronWise}. Elements in the numerical kernel
preserve an ample numerical class and are polarized
automorphisms after taking a multiple. Since \(H^0(X,T_X)=0\), the kernel is finite.
The lifted action is again proper and cocompact. Agol's theorem gives
virtual compact specialness \cite{CAgol}. A torsion-free finite-index
subgroup acts freely on the contractible \((\rho(X)-1)\)-dimensional
domain with closed manifold quotient, giving the stated dimension.

Finally a strictly convex symmetric cone of dimension at least three
is Lorentz. A finite-index subgroup of its integral divider preserves
the Lorentz form because unimodularity removes the scalar factor. By Borel
density the invariant form is unique up to scale, hence rational up to
scale since the representation is integral. In at least five variables
Meyer's theorem gives a rational isotropic line. Its integral orthogonal
group is then a nonuniform lattice, whereas it contains the dividing
uniform lattice with finite index. This contradiction proves
\(\rho(X)\le4\).
\end{proof}

\subsection*{Cohomological dimension of the numerical group}
The following lemma supplies the group to which
Theorem~\ref{thm:picard-boundedness} applies our conjecture.

\begin{lemma}[Dimension of a strictly convex quotient]\label{lem:divided-dimension}
Let \(\Lambda\) be a lattice of rank \(n\geq4\). Suppose
\(G\leq\mathrm{GL}(\Lambda)\) divides a full-dimensional properly
convex open domain \(U\subset\mathbb P(\Lambda\otimes\R)\) whose
boundary contains no nontrivial segment. There is a finite-index
subgroup \(G_0\leq G\cap\mathrm{SL}(\Lambda)\) that is finitely
generated, torsion-free and word-hyperbolic. Its given inclusion into
\(\mathrm{SL}(\Lambda\otimes\R)\) is projective Anosov and
\[
                          \operatorname{cd}_{\mathbb Z}(G_0)=n-1.
\]
\end{lemma}
\begin{proof}
Intersect \(G\) with the kernel of reduction modulo three in
\(\mathrm{SL}(\Lambda)\). This kernel is torsion-free and of finite
index in \(\mathrm{GL}(\Lambda)\). The resulting group \(G_0\) still
divides \(U\). Properness gives finite point stabilizers, so the
torsion-free action is free. Strict convexity of the divided domain
gives word hyperbolicity \cite[Theorem~3]{BenoistDiv}. The dividing
representation is projective Anosov
\cite[discussion before Corollary~1.4]{CanaryTsouvalas}.
In an affine chart \(U\) is convex and hence contractible. Thus
\(U/G_0\) is a closed aspherical \((n-1)\)-manifold with fundamental
group \(G_0\). A finite CW structure gives finite generation and the
upper bound \(n-1\) for cohomological dimension. Its top cohomology
with the orientation module is nonzero, which gives the reverse
inequality.
\end{proof}

\begin{corollary}[Equivalent projective formulation]\label{cor:projective-conjecture}
The Integral Anosov Conjecture is equivalent to the existence of an
absolute upper bound on the ranks \(n\) of lattices for which an
integral subgroup divides a full-dimensional properly convex
domain in \(\mathbb{RP}^{n-1}\) with strictly convex boundary.
\end{corollary}
\begin{proof}
If the conjecture holds, the lemma excludes such domains for
\(n>\max\{3,N\}\). Conversely, take \(n\geq4\) above the bound
on divided domains. A group satisfying the hypotheses of the conjecture
has cohomological dimension at most \(n-1\). If equality held,
Canary--Tsouvalas would give a strictly convex domain divided by
that very inclusion \cite[Theorems~1.3 and~1.5]{CanaryTsouvalas}.
This contradicts the bound on \(n\).
\end{proof}

The conjecture appears open in general. In the quadratic case,
a divided quadratic cone has a rational
defining form, up to scale, which is anisotropic over \(\mathbb Q\)
\cite[Theorem~2]{Tsuchihashi}. An indefinite rational form of rank at
least five is isotropic by Meyer's theorem \cite[Ch.~IV, \S3]{SerreCourse}.
Hence quadratic cones give \(n\leq4\). Choi--Choi's integral example in
\(\mathbb{RP}^4\) comes from a compact hyperbolic simplex group
\cite[Proposition~6.6 of the preprint]{ChoiChoi}. This gives equality
at \(n=5\), so any proposed \(N\) must be at least five.
Audibert--Douba--Lee--Marquis give a closed hyperbolic three-manifold
whose group virtually embeds Zariski-densely in
\(\mathrm{SL}(n,\mathbb Z)\) for every \(n\geq4\)
\cite[Theorem~1.4]{AudibertEtAl}. Its cohomological dimension stays
three as \(n\) grows. A family contradicting our conjecture would
have to realize \(n-1\) in unbounded matrix rank.

\begin{corollary}[Cocompact action on a convex subset]\label{cor:c-core}
Let a subgroup of \(\Gamma\) act cocompactly on a nonempty, invariant,
relatively closed convex subset \(K\subset\Omega\). If the subgroup
is not word-hyperbolic, \(X\) has a finite \'etale cover by an abelian
threefold. In particular \(\rho(X)\le9\).
\end{corollary}

\begin{proof}[Proof sketch]
The image of \(K\) in the quotient of the projectivized big cone by
this subgroup is compact and closed. Its ideal boundary therefore has
zero volume, since an accumulation point in the big cone would be in \(K\).
By Islam--Zimmer \cite[Theorem~5.1]{CIslamZimmer}, failure of hyperbolicity
gives a properly embedded simplex of dimension at least two. It contains
an open triangle whose three edges lie in the ideal boundary (partition
the vertices into three groups and take the three barycentres).
The triangle argument in Proposition~\ref{prop:c-group} therefore still
gives volume \(b t_1t_2t_3\) and places the triangle in one nef chamber
of a model \(Y\). On \(K\), compactness bounds
\(s/\operatorname{vol}^{1/3}\). Approaching each vertex within the
triangle gives \(k_Y(v_i)=0\). Their sum belongs to \(C\), so
nonnegativity of \(k_Y\) on \(M\) makes it identically zero. The finite
\'etale abelian cover then follows from the equality case for \(c_2\)
\cite{COguisoSakurai}. There are no rational curves on \(Y\), hence no
flops to another model. Pullback to the cover injects the divisor space
into a space of dimension at most \(h^{1,1}=9\).
\end{proof}

\begin{proposition}[A criterion for division of the movable domain]\label{prop:division-criterion}
Suppose no smooth projective flop model of \(X\) admits a divisorial
contraction or a fibration. Let \(\operatorname{Mov}^e(X)\) denote
the intersection of \(M\) with the cone of effective real divisor
classes. If \(\Gamma_X\Pi=\operatorname{Mov}^e(X)\) for a rational
polyhedral cone \(\Pi\subset\operatorname{Mov}^e(X)\), then
\(\Gamma_X\) divides \(\Omega_X\).
\end{proposition}
\begin{proof}
For any effective rational divisor \(D\) of nonzero numerical
class, run the minimal model program for \((X,\epsilon D)\), with
\(\epsilon>0\) small and rational. Each step is a flop: a first
divisorial step would occur on a smooth flop model. On the resulting
model the transform of \(D\) is nef and effective. Threefold log
abundance makes it semiample. Its nonzero class gives a positive-dimensional
image, and absence of a fibration forces this image to have dimension
three. Thus \(D\) is big and movable. This shows
\(M=\overline{\operatorname{Eff}}^{\,1}(X)\), so \(\operatorname{int}M\)
is the big cone. Each nonzero rational generator of \(\Pi\) has an
effective rational representative: a rational vector in the cone of
finitely many integral effective classes has a nonnegative rational
expression in those classes. All these generators therefore belong to
\(\operatorname{int}M\). Their convex hull in a positive affine slice
is compact and lies in \(\Omega_X\). Its translates cover \(\Omega_X\)
by hypothesis, proving cocompactness. Proper discontinuity follows
from discreteness of the integral projective group acting on a properly
convex domain. For the effective cone conjecture furnishing the cover,
see \cite{CGachet}.
\end{proof}

A divisorial contraction illustrates why the geometric hypothesis in
Theorem~\ref{thm:picard-boundedness}(2) still requires proof. For a divisorial contraction
\(f\) on a smooth model, the pullback of an ample divisor on its
three-dimensional target is big. It lies on the boundary of \(M\):
a curve sweeping the exceptional divisor has degree zero on that
pullback, positive degree on an ample class and nonnegative degree on
every movable class. For \(c_2(X)\) numerically nonzero, the equality
\(\operatorname{int}M=\BigCone(X)\) in the proof of
Proposition~\ref{prop:c-group} therefore precludes whole-domain
division. The argument above gives no uniform bound for nonfibered
infinite-chamber threefolds with such a contraction.

\section{Separating lines of plane cubic polars}\label{app:polars}

The coplanar case of Lemma~\ref{lem:div-walls} reduces to the following question about a real plane cubic. Can infinitely many of its polars have a line factor separating the parameter of that polar from all the others?

\begin{lemma}[Polars of a plane cubic]\label{lem:separating-polars}
Let \(F\) be a cubic form on a real three-dimensional vector space \(W\).
Suppose one of its polar quadratics \(B_e=\tfrac13\partial_eF\) is
nondegenerate. There is no infinite sequence of distinct directions
\([e_i]\in\mathbb P(W)\), no three collinear, with representatives
\(e_i\) and linear factors \(r_i\) of \(B_{e_i}\) satisfying
\[
                    r_i(e_i)<0<r_i(e_j)\qquad (i\ne j).
\]
\end{lemma}
\begin{proof}
Write \(T\) for the trilinear form of \(F\), so \(B_e(x)=T(e,x,x)\).
A conic with a line factor is singular. Thus all \([e_i]\) lie on the
Hessian curve \(\mathcal H=\{\det B_e=0\}\), a proper cubic since one
polar is nondegenerate.

First observe that a pencil through any point \([v]\) contains only
finitely many separating lines \(r_i=0\). Otherwise pass to the pencil
and descend its forms to \(W/\R v\). Fix \(r_1\) and normalize the
remaining vectors to \(r_1=1\), using positive scalars. Their images
on this affine line are distinct by separation. The middle of three
cannot be separated from the endpoints.

This observation already settles a reducible cubic \(F=\ell Q\) over
\(\R\). On \(\ker\ell\) we have \(3B_e=\ell(e)Q\). The binary
quadric \(Q|_{\ker\ell}\) is nonzero, since divisibility of \(F\) by
\(\ell^2\) would make it a cubic in at most two variables. Apart from
the at most two samples on \(\ell=0\), every factor line therefore
meets \(\ell=0\) in a zero of this fixed binary quadric. Its at most
two real zeros yield only finitely many separating lines.

Suppose now that \(F\) is irreducible. Pass to an accumulating
subsequence in an affine chart, rescaling all vectors positively after
fixing the sign of the chart. Write \(e_i=p+t_id_i\), where \(t_i>0\)
tends to zero and the unit secants \(d_i\) tend to a direction \(d\)
independent of \(p\). Normalize the factors by positive scalars so
\(r_i\to r\ne0\). For fixed \(i\), separation and the limit over
\(j\ne i\) give \(r_i(p)\ge0\) and \(r(e_i)\ge0\). Since
\(r_i(e_i)<0\), we obtain \(r(p)=0\) and
\[
                  r_i(d_i)<0\le r(d_i),\qquad r(d)=0.
\]
Also \(r\) divides \(B_p\), by compactness of projective factor lines.
Here \(B_p\ne0\), since otherwise symmetry would put \(p\) in the
kernel of every polar. Hence \(F(p)=0\). The line \(r=0\) through
\([p]\) contains the limiting secant to \(\mathcal H\).

If \([p]\) is smooth on \(F=0\), then \(B_p\) has rank two, since a square
vanishing at \(p\) would give \(T(p,p,-)=0\). Its kernel vector
\(k\) is independent of \(p\). The unique factor line through \([p]\)
is \(\mathbb P\langle p,k\rangle\), the tangent to \(F=0\).
On this line \(F(sp+tk)=t^3F(k)\), with \(F(k)\ne0\) by irreducibility.
The derivative of \(\det B\) in direction \(k\) is a nonzero multiple
of \(T(k,k,k)=F(k)\), by the adjugate formula for a rank-two matrix.
Thus \(\mathcal H\) is smooth at \([p]\) and transverse to this line,
contradicting the limiting secant in it.

If \([p]\) is singular on \(F=0\), every polar passes through \([p]\)
and \(B_p\) descends to a nonzero binary form on \(W/\R p\).
Remove the finitely many chosen lines passing through \([p]\) by the
pencil observation. For each remaining polar the complementary factor
passes through \([p]\). Whenever \(B_p(d_i,-)\ne0\), this factor is
exactly its tangent at \([p]\), since
\[
                  T(e_i,p,-)=t_i B_p(d_i,-).
\]
If \(B_p\) is a square up to scale, say \(B_p\propto\ell^2\), this
tangent is \(\ell\) for every sample off \(\ell=0\). Three such samples
span \(W\), so \(\ell\) divides every polar and hence \(F(x)=B_x(x)\).
If \(B_p=rs\) with distinct factors, then \(r(d)=0\) and \(s(d)\ne0\)
on the two-dimensional quotient. Hence the complementary factors tend
to \(B_p(d,-)\propto r\). The chosen factors also tend to \(r\), which
would make \(B_p\) a square. Both possibilities give a contradiction.
\end{proof}

\end{document}